\documentclass[]{article}

\input{epsf}
\usepackage{epsfig}
\usepackage{graphics}
\usepackage{amssymb,amsmath,amsthm,enumerate}
\usepackage{psfrag}
\usepackage{dcolumn}
\usepackage{color}
\usepackage{bm}
\usepackage{graphicx,amsmath,latexsym,enumerate,amsthm,color}
\usepackage{amssymb}
\usepackage{amsmath}
\usepackage[latin1]{inputenc}

\usepackage{pstricks,pst-plot,pst-text,pst-tree}
\usepackage{epic,eepic}
\usepackage{auto-pst-pdf}

\newcommand{\ti}{\tilde}

\newcommand{\RR}{I\!\!R}

\newtheorem{proposition}{Proposition}

\newtheorem{lemma}{Lemma}

\newcommand{\new}{}
\begin{document}

\begin{center}
{\bf \Large Gas-surface interaction and boundary conditions for the Boltzmann equation} \\
\vspace*{0.5 cm}

{\Large Stéphane Brull, Pierre Charrier, Luc Mieussens,\\

Institut Mathématiques de Bordeaux,\\
I.P.B.\\
Université de Bordeaux }\\

\today
\end{center}

{\bf Abstract}: In this paper we revisit the derivation of boundary
conditions for the Boltzmann Equation. The interaction between the
wall atoms and the gas molecules within a thin surface layer is
described by a kinetic equation introduced in~\cite{Borman} and used
in~\cite{ACD}. This equation includes a Vlasov term and a linear
molecule-phonon collision term and is coupled with the Boltzmann
equation describing the evolution of the gas in the bulk
flow. Boundary conditions are formally derived from this model by
using classical tools of kinetic theory such as scaling and systematic
asymptotic expansion. In a first step this method is applied to the
simplified case of a flat wall. Then it is extented to walls with
nanoscale roughness allowing to obtain more complex scattering
patterns related to the morphology of the wall. It is proved that the
obtained scattering kernels satisfy the classical imposed properties
of non-negativeness, normalization and reciprocity introduced by
Cercignani~\cite{C}.

\section{Introduction}
The Boltzmann equation is a powerful tool to describe phenomena in a gas flow  taking place
at a the scale of the order of the mean free path, i.e.  the
micrometric scale (for the air under stantard conditions). For many applications the gas flow takes place in a region bounded by
one or several solid bodies.  Then boundary conditions have to be prescribed in order to characterize the behavior of
the gas close to the wall~\cite{C,Sone2}.  \\
%This behavior is  \\

The first attempt to propose boundary conditions for the Boltzmann equation goes backs to
Maxwell in a paper of 1879 (\cite{Maxw}) where he discusses the way to
describe the interaction between a gas and a wall.
% Consider that the
%present gas is a bounded domain $\Omega$, whose boundary is $\partial
%\Omega$ and $\nu $ is the unit vector to the surface at $x$.
 The first condition he proposed corresponds to  a simple gas-solid interaction where we assume that the 
wall  is smooth, and perfectly elastic, so that the particles of gas are specularly
reflected.  This condition writes 
% the boundary conditions write
\begin{eqnarray}
\label{specf} f(t,x,v) = f(t,x,v-2 \nu \langle \nu ,v \rangle )
\hspace*{3 mm} \langle v , \nu \rangle >0,  
 \end{eqnarray} 
 where $\nu $ is the unit vector to the surface at point $x$ and
 $f(t,x,v)$ is the distribution function of particles that tt time $t$
 and position $x$ have the velocity $v$. Maxwell noticed that this assumption means that the gas can exert any
stress on the surface only in the direction of the normal. But this is
not physically relevant because in practical situations it can also
exert stress in oblique directions. 
This is why he introduced another type of boundary
conditions corresponding to a more complex gas-solid interaction. Physically he supposed that the wall has a stratum in
which fixed elastic spheres are placed. Moreover the stratum is
assumed to be deep enough so that every molecule going from the gas to
the wall must collide ones or more with the spheres.    
In this case, the particle is reflected into the gas with
a velocity taken with a probability whose density corresponds to the
equilibrium state of the gas. In that case the boundary condition (known as the perfect accommodation or diffuse reflexion condition) writes
\begin{eqnarray}
\label{maxdiff} f(t,x,v) =  \frac{1}{2 \pi (R T)} \int_{\langle v^{\prime}, \nu
  \rangle <0 } | \langle v^{\prime} , \nu  \rangle   |  \, f(t,x, v^{\prime}  )  dv^{\prime} \,  \exp (-\frac{v^2}{2RT}), 
\hspace*{3 mm} \langle v , \nu \rangle >0,   
 \end{eqnarray} 
where $T$ is the temperature of the wall. Finally Maxwell considered a more
complicated intermediate situation which is devoted to be more
physically realistic. This model is intermediate between the two
previous ones. Maxwell postulated that there is a fraction
of the gas which accomdates to the temperature of the solid and
another one which
is reflected by the solid. In that case the boundary conditions writes
\begin{eqnarray}
\nonumber f(t,x,v) &=& (1- \alpha ) f(t,x,v-2 \nu \langle \nu ,v \rangle )
\\ \label{maxwellbord} &+& \alpha \, \frac{1}{2 \pi (R T)} \int_{\langle v^{\prime}, \nu
  \rangle <0 } | \langle v^{\prime} , \nu  \rangle   |  \, f(t,x, v^{\prime}  )  dv^{\prime} \,  \exp (-\frac{v^2}{2RT}), 
\hspace*{3 mm} \langle v , \nu \rangle >0,   
 \end{eqnarray} 
where $T$ is still the temperature of the wall and $\alpha \in [0,1]$ is
called the accomodation coefficient. It represents the tendency of a gas to accomodate
to the wall. It means that a
fraction of $(1-\alpha)$ of molecules satisfies specular boundary
conditions whereas a fraction of $\alpha$ satisfies Maxwell diffuse boundary
conditions. When $\alpha = 0$, we recover the specular boundary
conditions and when $\alpha=1$, we recover the diffuse
boundary condition. The main drawback of this condition is that it gives  the same accommodation coefficient 
for energy and momentum though it is known that energy and momentum accommodate differently in physical 
molecule-wall interactions (see for instance
\cite{Cerlamp}). Nevertheless, this condition has been widely used,
both for theoretical studies and numerical simulations for practical applications.  \\

More recently, in \cite{CIL,C,Cersurf,Cer-inria} Cercignani adressed in great details the question of gas-surface 
interaction and boundary conditions for the Boltzmann equation with a large bibliography. He 
introduced a general formulation of the boundary conditions
\begin{equation}\label{eq-galbc} 
f(t,x,v) |\langle \nu ,v \rangle|_{|\langle \nu ,v \rangle >0}=\int_{\langle \nu ,v' \rangle<0} R(v^{\prime} \rightarrow v, x, t)f(t,x,v')|\langle \nu ,v' \rangle|dv',
\end{equation}
where the scattering kernel $R(v^{\prime} \rightarrow v, x, t)$
characterizes the interaction between the molecules of the gas and the
molecules of the wall. More precisely $R(v^{\prime} \rightarrow v, x,
t)$ represents the probability density that a molecule stricking the
wall with a velocity  $v^{\prime}$ at point $x$ and time $t$ is reemitted at the same point
with a velocity between $v$ and $v \, + dv$. To determine the
scattering kernel, Cercignani proposed
to use either physical or mathematical considerations. \\
In the physical approach we have to compute as exactly as possible the
path of the molecules within the wall.  This is anything but easy
since such a molecule may experience various events such as elastic
scattering, inelastic scattering (including multi-phonon scattering),
temporary or permanent adsorption, mobile adsorption (surface
diffusion), condensation, reactive interactions.  Therefore, in a
first attempt, very simplified models have been used to describe the
wall and the interactions such as arrays of smooth hard sphere or hard
cubes (see the work of Maxwell and the references given in
\cite{Cer-inria}). A more interesting way to approximate the path of
molecules within the wall has been proposed by Cercignani.  He
suggested to use a transport equation for the molecules inside the
solid which is regarded as a half-space.  This transport equation
includes a Vlasov-type term describing the van der Walls forces
exerted on the gas molecules by the solid atoms and a linear collision
term (of Boltzmann or Fokker-Plank type) describing the scattering
by phonons. Nevertheless, the Maxwell condition (\ref{maxwellbord})
can be recovered in this way (with a Boltzmann-like collision term) as
well as the Cercignani-Lampis condition \cite{Cerlamp} (with a
Fokker-Plank collision term). This latter condition is free from the
physically inconsistence of the Maxwell condition indicated above and
has been widely used.  More recent works come close to the same
approach by determining the molecule-wall interactions by means of
molecular dynamics simulation \cite{celestini,arya}. But an intrinsic
difficulty in this physical approach is due to our lack of knowledge
of the surface layers of solid walls, which leads Cercignani to
propose as an alternative that he called the mathematical approach.  \\
The idea of the mathematical approach is to construct a scattering
kernel, as simple as possible, satisfying the following basic
(physical)
requirements: \\
 
 \noindent
(i)  Non-negativeness:
 \begin{equation}
R(v^{\prime} \rightarrow v, x,t) \geq 0,  \label{nonneg}
\end{equation}
\noindent
(ii) Normalization: 
\begin{equation} \label{norm}
\int_{\langle v , \nu \rangle >0} R(v^{\prime} \rightarrow v, x,t) \,
dv =1 ,
\end{equation}
this property means  that the
mass flux through the boundary vanishes. It is valid when permanent adsorption is excluded.

\noindent
(iii) Reciprocity:  
\begin{eqnarray}
| \langle v^{\prime} , \nu \rangle | \, M_{w}(v^{\prime})  \,   
R(v^{\prime} \rightarrow v , x, t )
= | \langle v , \nu \rangle | \, M_{w}(v)  \,   
R(-v \rightarrow -v^{\prime} , x, t ), \label{recip}
\end{eqnarray}
where $ M_{w}$ is a Maxwellian distribution having the temperature of
the wall.  This last property means that the microscopic dynamics is
time reversible, and that the wall is in a local
equilibrium state and is not influenced by the incoming molecule. An example of a well-known scattering kernel derived  in such a way is the  Cercignani-Lampis model.  \\

In the present paper we use the so called physical approach but  we start from a somewhat more sophisticated  model 
introduced in \cite{Borman} and used in \cite{BKP2,Kry3,BBK,Kry1,Kry2,Kry4,ACD}
for studying gas-surface interaction, nanoflows and surface diffusion.  This model, valid for smooth walls,  is still a crude approximation of the complex
gas surface interaction, but  it proved to be remarkably useful  to give new insight on these issues.  It couples the Boltzmann 
equation in the bulk flow with a kinetic model inside a very thin surface layer (with width typically less than a nanometer) where   
the van der Waals forces are taken into account. This model includes a Vlasov term to take into account 
the part of the interaction potential that depends on the frozen position of the atoms  of the solid wall (the long range interactions), 
and a Boltzmann like linear collision term between molecules and phonons to take into account the thermal  fluctuations 
of the atoms of the solid  (short range interactions). 
%This is a three phase model (the bulk flow, the surface layer, the solid wall).
\\
It contains several characteristic times: the characteristic time of the Boltzmann equation in the bulk flow, the
characteristic time of the kinetic model in the surface layer, the characteristic time of flight of a molecule through
the surface layer, the characteristic molecule-phonon relaxation time.  Then using classical tools of kinetic theory
such as scaling asymptotic analysis we can derive various models corresponding to different regimes according to the
relative value of the characteristic times. Thus in \cite{ACD} surface kinetic and surface diffusion models have been
derived from this three phase model: they describe mobile adsorption and can be interpreted as {\it non local} boundary
conditions.  In the present paper, using different scalings, we derive {\it local} boundary conditions from the same
basic three phase model.  First, a weak molecule-phonon interaction regime is considered. In that case the particles of
the gas quickly cross the surface layer and the classical specular boundary condition is obtained. Then a strong
molecule-phonon interaction is investigated. In this situation the particles of the gas slowly cross the surface layer
and are thermalized by the wall leading to Maxwell-diffuse boundary conditions. Finally, an intermediate interaction is
assumed, and we get a Maxwell-like boundary condition (\ref{maxwellbord}), but with a fraction of diffusely evaporated
molecules that depends on the velocity. Moreover, the relationship between this coefficient and the surface-molecule
interaction potential is formulated.  One of the interesting asset of this boundary condition is that it gives different
accommodation coefficients for energy and (normal and tangential) momentum, contrary to the original Maxwell condition.
Moreover it must be noted that mobile adsorption (see \cite{ACD}) as well as elastic or inelastic scattering are treated
within the same framework.  Finally this analysis is extended to a non-smooth wall with nanoscale roughness assumed
to be periodic in the directions parallel to the surface. This leads
to a scattering kernel with more complex reflexion patterns that
depend on the wall morphology.
%We obtain again a Maxwell-like boundary condition with a modified scattering pattern 
%which is able to show peaks near certain limiting angles. 
\\

This paper is organized as follows. Section 2 deals with the presentation of nanoscale kinetic models describing the
interaction between a wall and particles in a very simplifed configuration with a flat wall and simplified expression of
the potential. In section 3, the boundary conditions are derived under these assumptions by using asymptotic analysis.
In section~\ref{sec:roughness}, the same analysis is extended to the more realistic case of a wall with nanoscale
roughness and a general potential.  Section~\ref{sec:conclusion} is devoted to some comments on these results and to
concluding remarks.

\section{Nanoscale kinetic models for gas-surface interaction}

In this section we recall the nanoscale models describing a gas flow near a wall introduced in \cite{Borman} and \cite{ACD}.
In these models the interaction between the wall and the gas molecules through Van der Walls forces are taken into account in a thin 
{\em surface layer} (with thickness $L$ typically smaller than one nanometer). In all the following, for the sake of simplicity, we assume 
that the molecules move in a 2D half-plane
 \footnote{As indicated in \cite{ACD} we can assume that the molecules move in the 3D half-space $(x,y,z), z<0$, provided that $f$
 is interpreted  as the marginal distribution function obtained by integrating the original  distribution function with respect to $v_y$.}
and we consider the following configuration: the solid 
is occupying the half-space $z>L$, the gas phase is constituted by the gas molecules in the half-space $z<0$, outside of the range 
of the surface forces, and we consider separately the {\em surface layer} $0<z<L$, where the gas molecules move within the range of 
the surface potential.
The gas flow in this surface layer is modelled by the 
collisionless Boltzmann equation (the size of this layer is much smaller than the mean free path of the molecules), with a Vlasov 
term to take into account the part of the interaction potential that depends on the frozen position of the atoms 
of the solid wall (the long range interactions), and a collision term between molecules and phonons to take into account the thermal 
fluctuations of the atoms of the solid (short range interactions) (see \cite{Borman} for a  physical justification of this approach). 
Since in many applications the surface potential is an attractive-repulsive potential, some of the molecules in the surface layer have
a total energy which is too small to escape from the potential well
and are {\em trapped} in the surface layer. On the other hand some
molecules, called the {\em free} molecules, have enough energy to escape from the potential well and can leave the surface layer and go into
the bulk flow. \\

Both type of molecules (trapped and free) are taken into account in this approach and we give now more details on the model describing 
their motion  in the surface layer.\\

\subsection{The surface potential}
We assume that the wall is flat and we use a 
simplified interaction potential which writes
\begin{equation}
{\cal V}(x,z)=%U(x)+
W(z), \label{potU+W}
\end{equation}
where 
%$U$ is a smooth potential and 
$W$ is an attractive-repulsive potential, ie;
\begin{enumerate}
\item[(H1)]   \;\; $0 \leq W(z)$ ,
\item[(H2)] \;\; $\lim_{z \rightarrow L}W(z) = + \infty,$
\item[(H3)]  \;\;  the potential $W$  is repulsive (i.e. $W'(z)>0$) for 
$z_m \leq z < L$ and is attractive $(W'(z)<0$) for $0<z<z_m$, and we
set $W(z_m)=0$
\item[(H4)]  \;\;
The range of the surface forces is finite and thus, the potential satisfies
$W(z)=W_m$ for $z<0$.
\end{enumerate}
% with $W(z)=W_m$ for $z\geq L$. Then if $W_m\gg U_m$ the well depth of the tangential component of the potential.
This simplified potential allows to uncouple the parallel motion and the normal motion of gas molecules near the solid wall, 
which makes the mathematical developments much easier. Moreover,
though not physically realistic, this potential is sufficient to obtain 
accurate information on the behavior of the gas near the walls (see \cite{ACD} for more details). Extension to a more realistic interaction potential
is considered in section~\ref{sec:roughness}.\\

It is useful to introduce in the surface layer the following velocity
variable, called equivalent velocity:
\begin{equation}\label{eq-ez} 
e_z=sign (v_z)\sqrt{v_z^2+2W(z)/m},
\end{equation}
which is the velocity of a particle whose total energy $\frac{1}{2}mv_z^2+W(z)$ would be a
kinetic energy $\frac{1}{2}me_z^2$ only. We denote by $e=(v_x,e_z)$ the corresponding two
dimensional velocity.

It will be more convenient to describe the distribution function of gas molecules in the surface layer as a function
of $e$ rather than a function of $(v_x,v_z)$.\\

\begin{figure}
  \centering
 \input{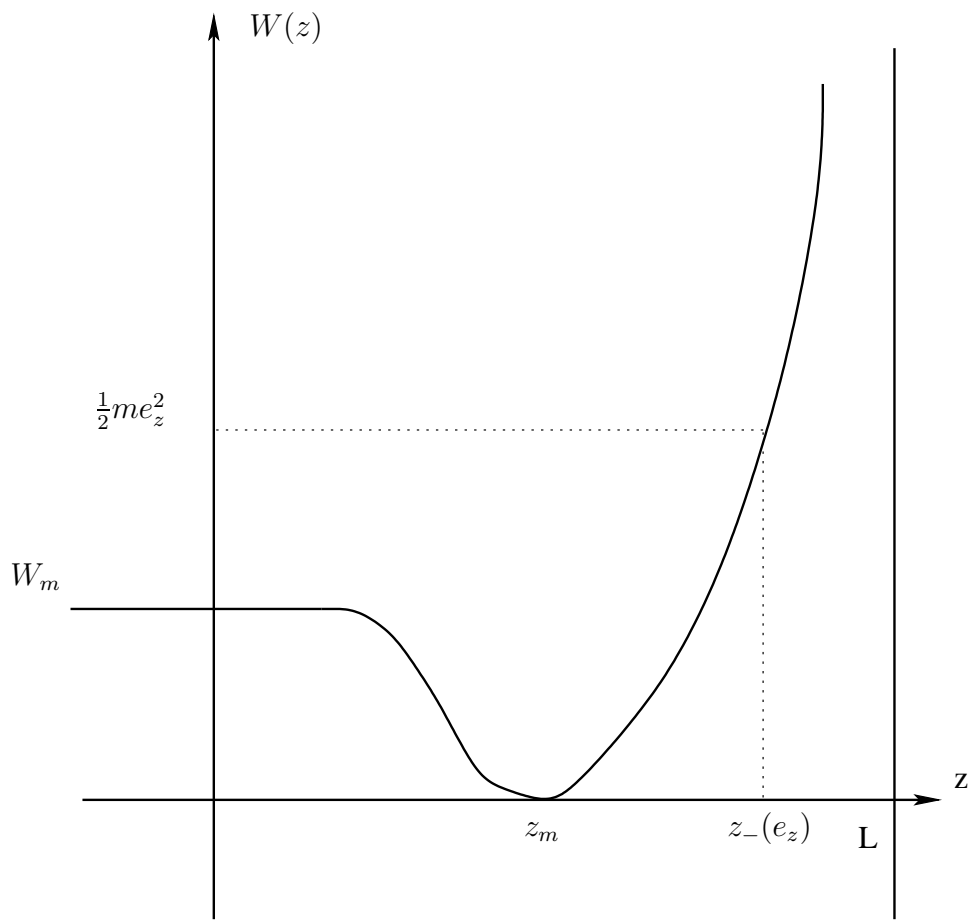}
\caption{Free particles: $|e_z|>\sqrt{\frac{2}{m}W_m}$ and $z<z_(e_z)$.
}
\label{figure:free}
\end{figure}

Now, we explain how particles can be divided into two different
classes: the free particles and the trapped particles. The
trajectory of a particle along $z$ is defined  (if there is no collision) by the two differential
equations $z'(t)=v_z(t)$ and $mv_z'(t)=\partial_zW(z(t))$. Along this
trajectory, the total energy $\frac{1}{2}mv_z^2+W(z)$ is
constant. According to the definition of the equivalent velocity
$e_z$ (see~(\ref{eq-ez})), we have
$\frac{1}{2}me_z^2=\frac{1}{2}mv_z^2+W(z)$ which is a constant too. A
particle is free if it can leave the surface layer and go into the
gas. In this case, the potential reaches the value $W_m$, and since
its kinetic energy $\frac{1}{2}mv_z^2$ is non-negative, this means
that $\frac{1}{2}me_z^2>W_m$, which is equivalent to
$|e_z|>\sqrt{\frac{2}{m}W_m}$. The limit position of this particle
when it is inside the surface layer is such that it takes a zero
velocity. At this point, denoted by $z_-(e_z)$, we have
$W(z_-(e(z)))=\frac{1}{2}me_z^2$ (see figure~\ref{figure:free}).

\begin{figure}
  \centering
 \input{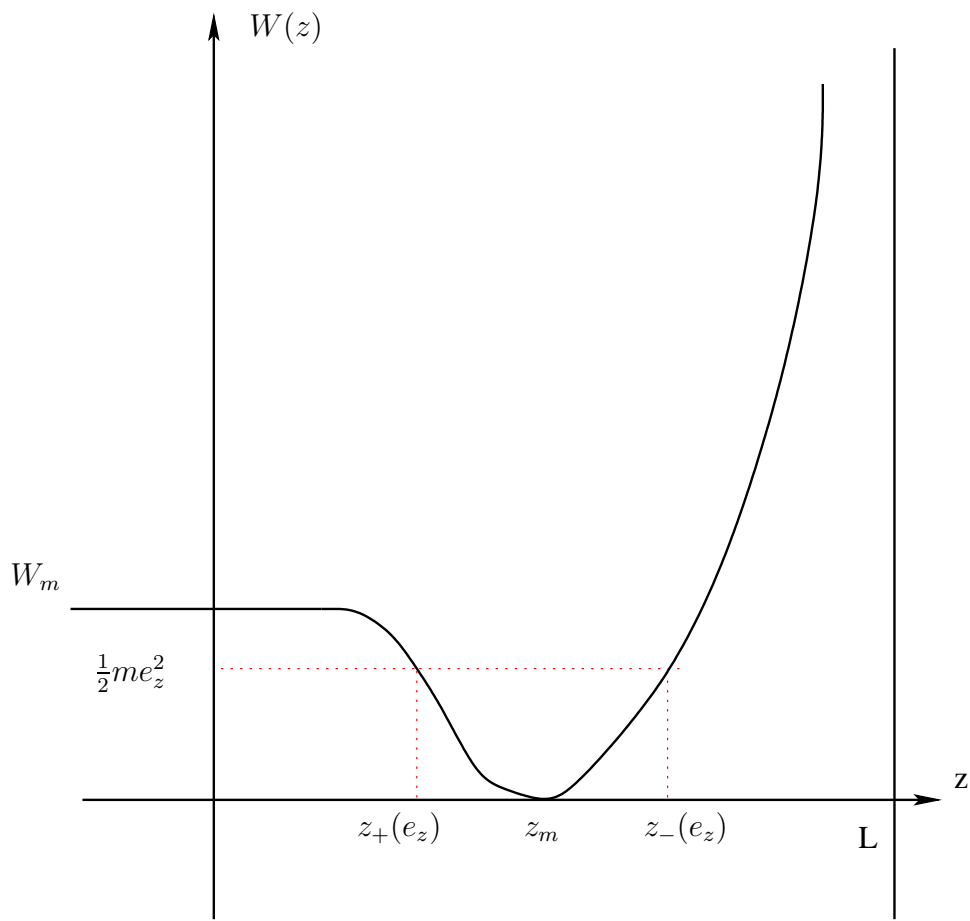}
\caption{Trapped particles: $|e_z|<\sqrt{\frac{2}{m}W_m}$ and $z_+(e_z)<z<z_-(e_z)$}
\label{figure:trapped}
\end{figure}

At the contrary, a particle is trapped if its total energy is lower
that $W_m$, that is to say $|e_z|<\sqrt{\frac{2}{m}W_m}$. In that
case, the potential is bounded by $\frac{1}{2}me_z^2<W_m$, which means
that $z$ varies between two limit values $z_+(e_z)$ and $z_-(e_z)$
such that $W(z_{\pm}(e_z))=\frac{1}{2}me_z^2$ (see
figure~\ref{figure:trapped}): the particle cannot escape from the
surface layer.

In order to have the same notation for trapped and free particles, we
set $z_+(e_z)=0$ for free particles (that is to say, if $|e_z|>\sqrt{\frac{2}{m}W_m}$). Moreover, for
particles with zero total energy, we have $e_z=0$ and hence the
velocity and the potential are zero too, which means that the particle
stay at position $z=z_m$. The we set $z_{\pm}(0)=z_m$ in this
case. 

With this definition, note that $z_+$ and $z_-$ are even functions of
$e_z$.\\

Now, we introduce some notations that are useful to switch between
$v_z$ and $e_z$ variables. The velocity of a particle with equivalent velocity $e_z$
located at position $z\in [z_{+}(e_z),z_{-}(e_z)]$ is 
given by
\begin{equation}
v_z(z,e_z)=sign (e_z)\sqrt{e_z^2-\frac{2}{m}W(z)},
\label{def-vz}
\end{equation}
and we have
\begin{equation}
v_z(z_{-}(e_z),e_z)=v_z(z_{-}(-e_z),-e_z)  =0. \label{cl-}
\end{equation}
Moreover, for trapped molecules we also have
\begin{equation}
v_z(z_{+}(e_z),e_z)=v_z(z_{+}(-e_z),-e_z) =0. \label{cl+}
\end{equation}
Let us define
\begin{equation*}
\sigma(z,e_z)= \frac{1}{|v_z(z,e_z)|} = (e_z^2-\frac{2}{m}W(z))^{-1/2}  \mbox{ for}\; |e_z| > \sqrt{2W(z)/m},\\
\end{equation*}
so that 
\begin{equation}  \label{eq-sigmavzez}
\sigma(z,e_z)\ v_z(z,e_z)= sign (e_z), 
\end{equation}
and also
\begin{eqnarray*}
{\tau}_z(e_z) &=& \int^{z_{-}(e_z)}_{z_{+}(e_z)} \sigma(z,e_z)dz
= \int^{z_{-}(e_z)}_{z_{+}(e_z)}(e_z^2-\frac{2}{m}W(z))^{-1/2}dz.
\end{eqnarray*}
As in \cite{DPV}, ${\tau}_z(e_z)$ can be interpreted as the time for a
molecule to cross the surface layer.  Moreover, for every $z \in
]0,L]$ the application $v_z \rightarrow e_z$ is a one-to-one
application from $[0,+\infty[$ onto $[ \sqrt{\frac{2}{m}W(z)}, +\infty
[$ and from $]-\infty, 0]$ onto $] -\infty,-
\sqrt{\frac{2}{m}W(z)}]$. Therefore differentiating (\ref{def-vz})
leads to
\begin{equation}\label{eq-chgvar1} 
dv_z = |e_z| \sigma(z,e_z)de_z.
\end{equation} 
 
\begin{figure}
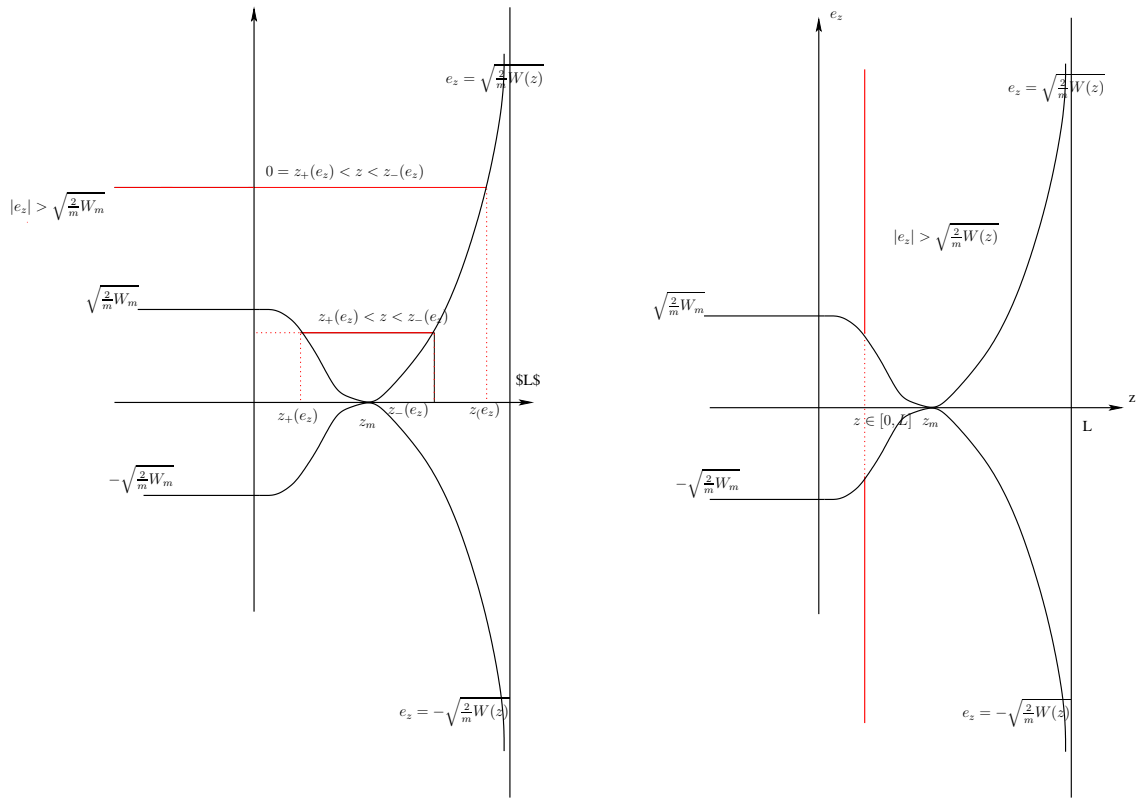

%  \centering
    \resizebox{!}{0.48\textheight} {\input{figurearticle1.pstex_t}} \hfill
    \resizebox{!}{0.48\textheight} {\input{figurearticle2.pstex_t}}
\caption{Domain of integration with respect to $e_z$ at a fixed $z$
  (left), with respect to $z$ at a fixed $e_z$ (right)}
\label{figure:fubini}
\end{figure}

Thus the integral of a given function $\psi(z,v_z)$ with respect to
$v_z$ can be transformed as follows:
\begin{equation}\label{eq-chgvar2} 
\int_{v_z} \psi(z,v_z)\ dv_z
= \int_{|e_z|>\sqrt{\frac{2}{m}W(z)}} \psi(z,v_z(z, e_z))\ |e_z| \sigma(z,e_z)de_z.
\end{equation}
Moreover, the order of integration in a $z$-$e_z$ integral can be
changed as follows (see figure~\ref{figure:fubini}): 
\begin{equation}\label{eq-fubini}
\begin{split}
& \int_0^L \left(\int_{|e_z|>\sqrt{\frac{2}{m}W(z)}} \psi(z,v_z(z, e_z))
  |e_z| \sigma_z(z,e_z)\ de_z \right)\ dz  \\
 & \qquad = \int_{-\infty}^{+\infty} \left( \int^{z_{-}(e_z)}_{z_{+}(e_z)}(\psi(z,v_z(z,
  e_z))|e_z| \sigma(z,e_z)\ dz \right)\ de_z. 
\end{split}
\end{equation}

\subsection{Molecule-phonon collision term}

In this paper we consider the general molecule-phonon collision term
\begin{equation*}
  Q[\phi](v)=\int_{\RR^2}K(v,v')\Bigl(\exp\Bigl({-\frac{m|v|^2}{2kT}}\Bigr)\phi(v')-\exp\Bigl({-\frac{m|v'|^2}{2kT}}\Bigr)\phi(v)\Bigr)\, dv'.
\end{equation*}
With the new velocity variable $e=(v_x,e_z)$ defined in~(\ref{eq-ez}),
for a given value of $z$, this operator reads:
\begin{equation}\label{eq-defQ} 
Q[\phi](z,e) = Q_+[\phi](z,e)-Q_-[\phi](z,e)=\int_{{\cal E}(z)} K(z,e,e')\left(G(e)\phi(e')- 
G(e')\phi(e)\right)J_{e'}\ de',
\end{equation}
where ${\cal E}(z)=\{ e',\ |e'_z|\geq\sqrt{2W(z)/m}  \},\; J_{e'}=J(z,e_z')=|e_z'|\sigma(z,e_z')$, and
\begin{equation}
G(e)=\exp \left(-\frac{m(|v_x|^2+|e_z|^2)}{2kT}\right).
\end{equation}
The collision kernel $K$ is such that $k(z,e \rightarrow e')=K(z,e,e')G(e')$
is the probability of transition per unit time from the state $e$ to the state $e'$  in a"collision" with a phonon.
The dimension of $K$ is $[time/length^2]$ (or, if the molecules move in a 3D plane, of $[time^2/length^3]$). 
We assume in the following 
%that $K$ does not depend on $z$ and
 that 
\begin{eqnarray}
\nonumber K(z,e,e') &=& K(z,e',e), \\
\label{Kbound} 0<\nu_0   &\leq&   K(z,e,e') \leq \nu_1, \\
K(z,v_x,-e_z,v'_x,-e'_z)&=&K(z,v_x,e_z,v'_x,e'_z), \label{Ksymez} \\
\nonumber K(z,-v_x,e_z,v'_x,e'_z)&=&K(z,v_x,e_z,v'_x,e'_z).
\end{eqnarray}

The loss term of the molecule-phonon collision term can be written
\begin{equation}\label{eq-lossterm} 
Q_-[\phi](z,e)=\frac{1}{\tau_{ms}(z,e)}\phi(z,e),
\end{equation}
where
\begin{equation}\label{eq-deftaums} 
\tau_{ms}(z,e)= \left( \int_{{\cal E}(z)}K(z,e,e')G(e')J({z,e_z'})de' \right)^{-1}
\end{equation}
is a collision time (at point $z$). It is useful for the sequel to introduce the mean relaxation time $\overline{\tau}_{ms}(e)$ defined as the harmonic mean of $\tau_{ms}(z,e)$ weighted
by $\sigma(z,e_z)$: 
\begin{equation}
\frac{1}{\overline{\tau}_{ms}(e)}=
\frac{\int_0^{z_-(e_z)}\sigma (z,e_z)/\tau_{ms}(z,e)\ dz}{\int_0^{z_-(e_z)}\sigma (z,e_z)\ dz}
=\frac{\int_0^{z_-(e_z)}\sigma (z,e_z)/\tau_{ms}(z,e)\ dz}{\tau_z(e_z)}. \label{meantaums}
\end{equation}
Using (\ref{Ksymez}), the even parity of $\sigma$, $G$, $J$ and
$z_{\pm}$ with respect to $e_z$, and the symmetry of ${\cal E}(z)$, we have :
\begin{equation*}
\tau_{ms}(z,v_x,-e_z)=\tau_{ms}(z, v_x,e_z), \mbox{ and}\; \; \overline{\tau}_{ms}(v_x,-e_z)=\overline{\tau}_{ms}( v_x,e_z).
\end{equation*}

Let us remark that if we assume $K(z,e,e')=1$, then $\tau_m$ does not
depends on $e$ and we have: 
\begin{equation}
Q[\phi]= \frac{1}{\tau_{ms}(z)}\left(\frac{n[\phi]}{\gamma(z)}G-\phi\right), \label{QBGK}
\end{equation}
where $\gamma (z)=\tau_{ms}(z)^{-1}=\int_{{\cal E}(z)}G(e')J(z,e_z')de'$ and $n[\phi]=\int_{{\cal E}(z)}\phi(e')J(z,e_z')de'$, which is quite similar to the BGK-like 
relaxation term used in \cite{ACD}. Finally we recall some of the main properties satisfied by the operator $Q$.

%%%%%%%%%%%%%%%%%%%%%%%%%%%%%%%%%%%%%%%%%%%%%%%%%%%%%
%%%%%PROPOSITION 1
%%%%%%%%%%%%%%%%%%%%%%%%%%%%%%%%%%%%%%%%%%%%%%%%%%%%
\begin{proposition}
The collision term satisfies the following properties
\begin{eqnarray}
\int_{{\cal E}(z)} Q[\phi](e)J_e\ de &=& 0,\; \mbox{ (mass conservation)},\label{eq-cons} \\
Q[\phi]=0 &\Leftrightarrow& \phi = n\ G, \; \mbox{(equilibrium)},\\
%Im\ Q &=& \{ \varphi ,\ n[\varphi]=\int_{{\cal E}(z)} \varphi (e)J_ede =0\},\\
 \int_{{\cal E}(z)}  Q[\phi](e)\phi(e)\frac{J_e}{G(e)}de  &\leq& 
-\nu_0\gamma(z)\int_{{\cal E}(z)} w^2 \frac{J_e}{G}de,\   \mbox{(H theorem)}, \\
\int_{{\cal E}(z)} Q[\phi](e)\psi(e)\frac{J_e}{G(e)}de &=& \int_{{\cal E}(z)}Q[\psi](e)\phi(e)\frac{J_e}{G(e)}de,\; \mbox{(symmetry) },
\end{eqnarray}
where we used the macro-micro decomposition $\phi=q+w$ with
$q=n[\phi]G$ and where $w=\phi-q$ satisfies $n[w]=0$.
\end{proposition}
\vspace*{0.5cm}

\subsection{Nanoscale models}
The first model introduced in \cite{Borman} and \cite{ACD} is the following system of coupled kinetic equations
which describes the flow of molecules in the surface layer (where the Van der Waals forces are acting) and outside:
\begin{eqnarray}
\partial_tf+v_x\partial_xf+v_z\partial_zf &=& 0,\ z<0 \label{boltzmannmax}  \\
f(t,x,0,v_x,v_z)|_{v_z<0}&=&\phi(t,x,0,v_x,e_z(0,v_z)),\label{interfacevzposmax}\\
%\int v_z f(t,x,0,v_x,v_z)\ dv&=& \int_{{\cal E}(z)} e_z \phi(t,x,0,v_x,e_z)\ de, \label{contfluxmax} \\
\partial_t\phi+v_x\partial_x\phi
%-\frac{U'(x)}{m}\partial_{v_x}\phi 
+ v_z(z,e_z)\partial_z\phi &=& Q[\phi]
,\ z_{+}(e_z)<z<z_{-}(e_z), \label{surfacelayer}\\ 
 \phi(t,x,0,v_x,e_z)_{e_z>\sqrt{2W_m/m}}&=&f(t,x,0,v_x,v_z(0,e_z)),\label{interfacevznegmax}\\
\phi(t,x,z_{-}(e_z),v_x,e_z)&=& \phi(t,z_{-}(-e_z),v_x,-e_z), \label{reflexion} \\
\phi(t,x,z_{+}(e_z),v_x,e_z)&=& \phi(t,z_{+}(-e_z),v_x,-e_z),\   |e_z|<\sqrt{2W_m/m},    \label{reflexion2}
\end{eqnarray}
where $f=f(t,x,z,v_x,v_z)$ is the distribution function describing the bulk flow and $\phi=\phi(t,x,z,v_x,e_z)$
is the distribution function describing the gas flow inside the surface layer. Let us remark that since we have chosen 
to define $\phi$ as a function of $(v_x,e_z)$ equation (\ref{surfacelayer}) does not contain a Vlasov term in 
the z-direction.\\

The above model describes the gas-solid interaction at the nanoscale, i.e on a domain $[0,x^*] \times [-z^*,L]$ 
with $x^*$ and $z^*$ $\approx$ 1 nanometer. But on a larger scale in the tangential direction, this model is too 
complicated and contains stiff terms that would make its numerical solution too much expensive. Thus in
\cite{ACD} the authors derived a limit model obtained by asymptotic
analysis when the domain is much larger than the surface layer (that
is to say $x^{\star}\approx z^* \gg L$).
In this model, the flow of molecules in the surface layer is described by a one-dimensional kinetic equation which can be 
considered as a nonlocal boundary condition for the Boltzmann equation in the bulk flow. \\

But on a larger scale in $x$ and $z$ this last model is still too complicated to manage and it would be interesting to 
investigate the relation between these nanoscale models and the standard boundary conditions used with 
the Boltzmann equation in gas kinetic theory. \\

In the following, we use the nanoscale model
(\ref{boltzmannmax}-\ref{reflexion2}) to derive various boundary conditions for the 
Boltzmann equation (\ref{boltzmannmax}), according to convenient scalings.

\section{Derivation of boundary conditions: case of a flat wall}
In this section, we assume that the characteristic times of the flow
in the surface layer (the time for a molecule to cross the surface
layer and the relaxation time of molecules by phonons) are much
smaller than the characteristic time of evolution of the bulk flow. We
derive boundary conditions for the Boltzmann equation in the bulk flow
by an asymptotic analysis of system
(\ref{boltzmannmax}-\ref{reflexion2}). The main point in this derivation is
to find the solution of a linear kinetic
problem which describes, in a first approximation, the motion of the molecules in the surface layer. Unfortunately this problem cannot be solved exactly but approximated solutions can be obtained (see Lemma 1) through an iterative process. 

We consider system (\ref{boltzmannmax}-\ref{reflexion2}) and we
introduce the following dimensionless quantities: 
\begin{eqnarray*}
\ti{n}=\frac{n}{n^*}, \ti{v}_{x}=\frac{v_{x}}{v^*}, \ti{v}_{z}=\frac{v_{z}}{v^*}, \ti{e}_{z}=\frac{e_{z}}{v^*}, \ti{f}=\frac{f}{f^*}, \ti{\phi}=\frac{\phi}{f^*},
 \ti{x}=\frac{x}{l^*},  \\
%\ti{U}=\frac{U}{U^*}, 
\ti{W}=\frac{W}{W^*}, \ti{W}_m=\frac{W_m}{W^*}, \ti{t}=\frac{t}{t_{B}^*},  \ti{\tau}_z=\frac{\tau_z}{\tau_z^*},  \ti{\tau}_{ms}=\frac{\tau_{ms}}{\tau_{ms}^*}, 
\ti{K}=\frac{K}{K^*},
\end{eqnarray*}
and $\ti{z}=\frac{z}{l^*}$ for the Boltzmann equation in the bulk
flow, while $\ti{z}=\frac{z}{L}$ in the surface layer. The reference
quantities are the followings: $n^*$ is the reference number density,
$v^*=\sqrt{kT/m}$, $f^*=n^*/{v^{*}}^2$, $t_B^*$ is the
reference time of evolution for the Boltzmann equation
(\ref{boltzmannmax}), $l^*=v^*t_B^*$, $\tau_{ms}^*=1/(K^*{v^*}^2)$ is a reference
relaxation time, $\tau_z^*=L/v^*$ is the characteristic time of
flight of a molecule through the surface layer,  and
$W^*=m{v^*}^2/2$.  

% We are now considering three different regimes
% according to the respective values of the three different
% characteristic times appearing in this problem, $\tau_z^*$,
% $\tau_{ms}^*$ and $t_B^*$. First, we consider the case where
% $\tau_{ms}^* =c_0 \ \tau_z^* \ll t_B^*$,
% where $c_0$ is a positive constant, and in a second step we consider the limiting regimes corresponding respectively to $c_0=0$ or $c_0=+\infty$. \\

In order to study different regimes corresponding to different order
of magnitude of the characteristic time scales $\tau_z^*$,
 $\tau_{ms}^*$ and $t_B^*$, we introduce the following nondimensional
 parameters:
\begin{equation*}
\varepsilon= \frac{\tau_{ms}^*}{t_B^*} \qquad \text{ and } \qquad \eta=\frac{\tau_{ms}^*}{\tau_z^*}. 
\end{equation*}
Then system (\ref{boltzmannmax}-\ref{reflexion2}) reads
in dimensionless form
\begin{eqnarray}
\partial_{\ti{t}}\ti{f}+\ti{v}_x\partial_{\ti{x}}\ti{f}+\ti{v}_z\partial_{\ti{z}}\ti{f} &=& 0,\ \ti{z}<0, \label{boltzmann-adim-2}  \\
\ti{f}(\ti{t},\ti{x},0,\ti{v}_x,\ti{v}_z)_{\ti{v}_z<0}&=&\ti{\phi}(\ti{t},\ti{x},0,\ti{v}_x,\ti{e}_z(0,\ti{v}_z)),
\label{interfacevzposadim-2}\\
%\int \int \ti{v}_z \ti{f}(\ti{t},\ti{x},0,\ti{v}_x,\ti{v}_z)\ d\ti{v}_xd\ti{v}_z&=& 
%\int_{|\ti{e_z}|>\sqrt{\ti{W}_m}} \int \ti{e}_z \ti{\phi}(\ti{t},\ti{x},0,\ti{v}_x,\ti{e}_z)\ d\ti{v}_xd\ti{e}_z , \label{contfluxadim-2} \\
\partial_{\ti{t}}\ti{\phi}+\ti{v}_x\partial_{\ti{x}}\ti{\phi}
%-\frac{\ti{U}'(\ti{x})}{2}\partial_{\ti{v}_z}\phi 
+ \frac{\eta}{\varepsilon}\ti{v}_z(\ti{z},\ti{e}_z)\partial_{\ti{z}}\ti{\phi} &=& 
\frac{1}{\varepsilon}\ti{Q}[\ti{\phi}],\ \ti{z}_{+}(\ti{e}_z)<\ti{z}<\ti{z}_{-}(\ti{e}_z), \label{reflexionzmadim}\\ 
 \ti{\phi}(\ti{t},\ti{x},0,\ti{v}_x,\ti{e}_z)_{\ti{v}_z>0}&=&\ti{f}(\ti{t},\ti{x},0,\ti{v}_x,\ti{v}_z(0,\ti{e}_z)),\label{interfacevznegadim-2}\\
\ti{\phi}(\ti{t},\ti{x},\ti{z}_{-}(\ti{e}_z),\ti{v}_x,\ti{e}_z)&=& \ti{\phi}(\ti{t},\ti{z}_{-}(\ti{e}_z),\ti{v}_x,-\ti{e}_z),\\
\ti{\phi}(\ti{t},\ti{x},\ti{z}_{+}(\ti{e}_z),\ti{v}_x,\ti{e}_z)&=& \ti{\phi}(\ti{t},\ti{z}_{+}(\ti{e}_z),\ti{v}_x,-\ti{e}_z),
\mbox{ for}\; |\ti{e}_z| < \sqrt{\ti{W}_m}.\label{reflexionzpadim}
\end{eqnarray}
We mention that with this dimensionless variables, a particle of
velocity $e_z$ located at $z$ is: 
\begin{itemize}
  \item either trapped if $|\ti{e}_{\ti{z}}|<\sqrt{\ti{W}_m}$, and hence stays between $\ti{z}_{\pm}(\ti{e}_{\ti{z}})$
defined by $\ti{W}(\ti{z}_{\pm}(\ti{e}_z))=\ti{e}_z^2$,
\item or free if $|\ti{e}_z|>\sqrt{\ti{W}_m}$, and hence stays on the
  left-hand-side of $\ti{z}_{-}(\ti{e}_z)$
defined by $\ti{W}(\ti{z}_{-}(\ti{e}_z))=\ti{e}_z^2$. We set
$\ti{z}_{+}(\ti{e}_{\ti{z}})=0$ in this case.
\end{itemize}

We can obtain boundary conditions for the Boltzmann equation through an asymptotic analysis of the above system when $\varepsilon \rightarrow 0$. This leads to the following results.

%%%%%%%%%%%
%%PROPOSITION 2
%%%%%%%%%%%
\begin{proposition}
Under the hypothesis (\ref{potU+W}) and (H1-H4), in the limit 
$\varepsilon\rightarrow 0$, 
the gas-surface interaction depends on the order of magnitude of
$\eta$ and can be described by the following boundary
conditions at $z=0$:
\begin{enumerate}
\item for $\eta=O(\frac{1}{\varepsilon})$, the boundary
  condition is 
%$limit $\frac{\tau_z^*}{\tau_{ms}^*}\rightarrow 0$ (or
%$c_0\rightarrow +\infty$), we recover the 
the specular reflection
\begin{eqnarray*}
{f}({t},{x},0,{v}_x,{v}_z)_{|{v}_z<0}&=&{f}({t},{x},0,{v}_x,-{v}_z).
\end{eqnarray*}
\item for $\eta=O(\varepsilon)$, the boundary condition is 
%In the limit $\frac{\tau_{ms}^*}{\tau_z^*}\rightarrow 0$ (or
%$c_0\rightarrow 0)$, we recover 
the reflection with perfect accommodation
\begin{eqnarray}
{f}({t},{x},0,{v}_x,{v}_z)_{|{v}_z<0}&=&{\kappa}(t,x) M(v_x,v_z), \label{perfaccomodationcond}
\end{eqnarray}
where
\begin{eqnarray*}
\kappa (t,x) &=& \int_{v_z>0}\int v_z{f}({t},{x},0,{v}_x,{v}_z)dv_xdv_z/
\int_{v_z>0}\int v_zM({v}_x,{v}_z)dv_xdv_z
\end{eqnarray*}
is such that the mass flux of $f$ through the boundary $z=0$ is zero,
and where $M(v)=\exp \left( -m(v_x^2+v_z^2)/2kT\right)$.\

\item for $\eta=O(1)$, the boundary condition writes, in a first approximation, as a  Maxwell-like boundary condition
\begin{equation*} \label{Maxcond} 
{f}({t},{x},0,{v}_x,{v}_z)_{|{v}_z<0} = {a}({v})\beta_1(t,x) M(v)+(1-{a}({v})){f}({t},{x},0,{v}_x,-{v}_z), 
\end{equation*}
and
% \begin{eqnarray}
% a(v) &=& 1-\exp \left( -\frac{2\tau_z(e_z(0,v_z))}{\overline{\tau}_{ms}(v_x,e_z(0,v_z))}\right)  \label{accomodation} \\
% \nonumber \beta_1(t,x) &=& \int_{v_z>0}\int v_za(v){f}({t},{x},0,{v}_x,{v}_z)dv_xdv_z/
% \int_{v_z>0}\int v_za(v)M({v}_x,{v}_z)dv_xdv_z,
% \end{eqnarray}
\begin{eqnarray}
a(v) &=& 1-\exp \left( -\frac{2\hat{\tau}_z(v_z)}{\hat{\tau}_{ms}(v)}\right)  \label{accomodation} \\
\nonumber \beta_1(t,x) &=& \int_{v_z>0}\int v_za(v){f}({t},{x},0,{v}_x,{v}_z)dv_xdv_z/
\int_{v_z>0}\int v_za(v)M({v}_x,{v}_z)dv_xdv_z,
\end{eqnarray}
with the notations $\hat{\tau}_z(v_z)=\tau_z(e_z(0,v_z))$ and
$\hat{\tau}_{ms}(v)=\overline{\tau}_{ms}(v_x,e_z(0,v_z))$. 
This boundary condition ensures a zero mass flux of $f$ at the
boundary $z=0$. Moreover, it can be written under the
general form~(\ref{eq-galbc}) with a scattering kernel $R(v'\rightarrow v)$ that
satisfies the properties of non-negativeness, normalization and reciprocity.
\end{enumerate}

\end{proposition}

{\em Proof:} \\
In order to simplify the notations, the tilde $\tilde{}$ over the
dimensionless quantities are dropped in the following. To
 avoid confusion, we will indicate explicitely when we come back to dimensional quantities.\\

In order to perform an asymptotic analysis of system
(\ref{boltzmann-adim-2}-\ref {reflexionzpadim}), we look for a
solution in the form
\[{f}={f}_{\varepsilon}={f}^0+\varepsilon {f}^1+ ...,\;\;\;
{\phi}={\phi}_{\varepsilon}={\phi}^0+\varepsilon {\phi}^1+ ...  . \]
This expansion is inserted
into~(\ref{boltzmann-adim-2}--\ref{reflexionzpadim}) and we
identify the terms of same power of magnitude w.r.t $\varepsilon.$ The
zeroth-order term $f^0$ satisfies 
\begin{eqnarray}
\partial_{{t}}{f}^0+{v}_x\partial_{{x}}{f}^0+{v}_z\partial_{{z}}{f}^0 &=& 0 ,\label{eq-f0} \\
{f}^0({t},{x},0,{v}_x,{v}_z)_{{v}_z<0}&=&{\phi}^0({t},{x},0,{v}_x,{e}_z(0,v_z)).
\label{orderzerointerfacevzpos}
\end{eqnarray}
However, the zeroth-order term $\phi^0$ depends on the order of
magnitude of $\eta$. \\

(1) We consider the case
$\eta=O(\frac{1}{\varepsilon})$. This means that
 \[  \tau_z^* \ll \tau_{ms}^* \ll t_B^* , \] that is to say the free time of flight of a molecule 
to cross the surface layer is much smaller than the relaxation time of
molecules by phonons. Thus the flow of molecules
crosses the surface layer so quickly that the relaxation phenomena can be neglected.
Then $\phi^0$ satisfies the following linear kinetic surface layer (LKSL) problem:
 \begin{eqnarray}
  {v}_z({z},{e}_z)\partial_{{z}}{\phi}^0 &=& 0,\; 
 \mbox{ for}\ {z}_+({e}_z)<{z}<{z}_-({e}_z),\  \label{spec}\\
{\phi}^0({t},{x},0,{v}_x,{e}_z)_{{e}_z>\sqrt{{W}_m}}&=&{f}^0({t},{x},0,{v}_x,{v}_z(0,{e}_z)),
\label{interfacevznegspec}\\
{\phi}^0({t},{x},{z}_{-}({e}_z),{v}_x,{e}_z)&=& {\phi}^0({t},{x},{z}_{-}({e}_z),{v}_x,-{e}_z),
\label{orderzeroreflexionzmspec}\\
{\phi}^0({t},{x},{z}_{+}({e}_z),{v}_x,{e}_z)&=& {\phi}^0({t},{x},{z}_{+}({e}_z),{v}_x,-{e}_z), 
\mbox{ for}\; |{e}_z| < \sqrt{{W}_m}. \label{orderzeroreflexionzpspec}
\end{eqnarray}
Consider some $v_z<0$ and the boundary
condition~(\ref{orderzerointerfacevzpos}) where we write
$\bar{e}_z={e}_z(0,{v}_z)$: 
\begin{equation}\label{eq-CLfb}
{f}^0({t},{x},0,{v}_x,{v}_z)_{v_z<0} = {\phi}^{0}({t},{x},0,{v}_x,\bar{e}_z).
\end{equation}
Since~(\ref{spec}) implies that
${\phi}^{0}$ does not depend on $z$, we can replace $z=0$ in the
right-hand side of~(\ref{eq-CLfb}) by $z={z}_{-}({e}_z)$ to get
\begin{eqnarray*}
{f}^0({t},{x},0,{v}_x,{v}_z)_{v_z<0} = {\phi}^{0}({t},{x},{z}_{-}(\bar{e}_z),{v}_x,\bar{e}_z).
\end{eqnarray*}
Moreover~(\ref{orderzeroreflexionzmspec}) and the even parity of $z_-$ imply
\begin{eqnarray*}
{f}^0({t},{x},0,{v}_x,{v}_z)_{v_z<0} = {\phi}^{0}({t},{x},{z}_{-}(-\bar{e}_z),{v}_x,-\bar{e}_z).
\end{eqnarray*}
Again, we use the fact that ${\phi}^{0}$ does not depend on $z$ to get
\begin{eqnarray*}
{f}^0({t},{x},0,{v}_x,{v}_z)_{v_z<0} = {\phi}^{0}({t},{x},0,{v}_x,-\bar{e}_z),
\end{eqnarray*}
where, by definition, $-\bar{e}_z\geq \sqrt{W_m}$.  Now we can
use~(\ref{interfacevznegspec}) to replace the right-hand side of the
previous relation and to get the specular boundary condition 
\begin{eqnarray}
 {f}^0({t},{x},0,{v}_x,{v}_z)_{v_z<0}	
&=&   {f}^0({t},{x},0,{v}_x,-{v}_z).\;           \nonumber       
\end{eqnarray}
We mention that we used $\eta=O(\frac{1}{\varepsilon})$ for
simplicity. In fact, we recover the same boundary condition if
$\eta=O(\varepsilon^{-\alpha})$, for every positive $\alpha$.

(2) Now we assume $\eta=O(\varepsilon)$, which implies that
\[ \tau_{ms}^*\ll \tau_z^*  \ll t_B^*.\]
This means that the relaxation time of molecules by phonons is much smaller than 
the free time of flight of a molecule to cross the surface layer. In this limit the flow of incoming molecules into the 
surface layer immediately relaxes toward the equilibrium. 
% We introduce the small parameter
% \begin{equation*}
% \varepsilon= \frac{\tau_z^*}{\tau_{ms}^*}=\frac{\tau_{ms}^*}{t_B^*}.
% \end{equation*}
Now the LKSL problem satisfied by $\phi^0$ reads
\begin{eqnarray}
\  {Q}[{\phi}^0] &=& 0, \; \; \mbox{ for}\   {z}_+({e}_z)<{z}<{z}_-({e}_z),
\label{orderzeroperf}\\
{\phi}^0({t},{x},0,{v}_x,{e}_z)_{{e}_z>0}&=&{f}^0({t},{x},0,{v}_x,{v}_z(0,{e}_z)),
\label{orderzerointerfacevznegperf}\\
{\phi}^0({t},{x},{z}_{-}({e}_z),{v}_x,{e}_z)&=& {\phi}^0({t},{z}_{-}({e}_z),{v}_x,-{e}_z),
\label{orderzeroreflexionperf} \\
{\phi}^0({t},{x},{z}_{+}({e}_z),{v}_x,{e}_z)&=& {\phi}^0({t},{x},{z}_{+}({e}_z),{v}_x,-{e}_z), 
\mbox{ for}\; |{e}_z| < \sqrt{{W}_m} \label{orderzeroreflexionzpperf}
\end{eqnarray}
which gives
\begin{equation}\label{eq-phi0maxw} 
{\phi}^{0}({t},{x},{z},{v}_x,{e}_z)=\alpha (t,x) G(v_x,e_z),\, \mbox{ for}\  {z}_+({e}_z)<{z}<{z}_-({e}_z).
\end{equation}
However, the distribution function ${\phi}^0$ is Maxwellian and hence
cannot satisfy the inflow boundary condition (\ref{orderzerointerfacevznegperf}). Thus we have to introduce in  
the expansion of ${\phi}$ a Knudsen-layer corrector
\begin{equation*}
{\phi}({t},{x},{z},{v}_x,{e}_z)={\phi}^0({t},{x},{z},{v}_x,{e}_z)+
{\psi}^0({t},{x},\frac{{z}}{\varepsilon},{v}_x,{e}_z)+
\varepsilon{\phi}^1({t},{x},{z},{v}_x,{e}_z)+ ...,
\end{equation*}
where ${\phi}^0$ is still defined by~(\ref{eq-phi0maxw}) and
satisfies~(\ref{orderzeroperf}, \ref{orderzeroreflexionperf},
\ref{orderzeroreflexionzpperf}), and ${\psi}^0({t},{x},y,{v}_x,{e}_z)$ is given by
\begin{eqnarray}
{v}_z(0,{e}_z)\partial_y {\psi}^0 &=& {Q}[{\psi}^0],\; \mbox{ for } |{e}_z|\geq \sqrt{{W}_m},\
0<y<+\infty, \label{Kcorr1}\\
{\psi}^0({t},{x},0,{v}_x,{e}_z)_{|{e}_z>0} &=& {f}^0({t},{x},0,{v}_x,{v}_z(0,{e}_z))-
{\phi}^0({t},{x},0,{v}_x,{e}_z),\label{Kcorr2} 
%\\
%{\psi}^0({t},{x},+\infty,{v}_x,{e}_z) &=& 0 \label{Kcorr3}
\end{eqnarray}
and should rapidly decrease to $0$ for large $y$. 
% Let us note that 
% \begin{equation}
% \int \int_{{\cal E}_z({z})}{e}_z{\psi}^0({t},{x},0,{v}_x,{e}_z)\ d{v}_xd{e}_z =0.
% \label{massflux2}
% \end{equation}
It is useful to introduce ${\chi}({t},{x},y,{v}_x,{e}_z)$
defined by
\begin{equation}\label{eq-defchi} 
{\chi}({t},{x},y,{v}_x,{e}_z)= {\psi}^0({t},{x},y,{v}_x,{e}_z)
+\phi^0({t},{x},0,{v}_x,{e}_z).
\end{equation}
% so that  
% \begin{equation}
% {f}^0({t},{x},0,{v}_x,{v}_z)_{|{v}_z<0}={\chi}({t},{x},0,\ti{v}_x,{e}_z(0,{v}_z)). 
% \label{f0out}
% \end{equation}
 Thus ${\chi}$ is the unique bounded solution of the following linear half-space problem
\begin{eqnarray}
{v}_z(0,{e}_z)\partial_y {\chi} &=& {Q}_0[{\chi}]
%={Q}_{+,0}[{\chi}]-\frac{{\chi}}{{\tau}_{ms}(0,{e})},
\label{chi1} \\
{\chi}({t},{x},0,{v}_x,{e}_z)_{|{e}_z>0} &=& {f}^0({t},{x},0,{v}_x,{v}_z(0,{e}_z)), \label{chi2} 
\end{eqnarray}
where
\begin{equation*}
{Q}_0[{\chi}] = \int_{{\cal E}(0)} {K}(0,{e},{e}')\left({G}({e}){\chi}({e}')- 
{G}({e}'){\chi}({e})\right){J}(0 ,{e}')\ d{e}'.
\end{equation*}
Then using the result given in \cite{Golse-Klar} and \cite{Klar}
on the linear half-space problem we get the approximation
\begin{equation}
{\chi}({t},{x},y,{v}_x,{e}_z)_{|{e}_z<0} \approx {\chi}^{(1)}({t},{x},y,{v}_x,{e}_z)_{|{e}_z<0}
= \kappa(t,x){G}\label{approxsol}
\end{equation}
for every $y\geq 0$, where $\kappa$ can be determined as follows. A
standard result on the linear half-space
problem~(\ref{chi1}--\ref{chi2}) shows that 
$\chi$ necessarily satisfies $\int_{{\cal E}(z)}e_z
\chi(t,x,y,v_x,e_z) \, de=0$ for every $y$. Then, writing this relation at
$y=0$ and using the boundary condition~(\ref{chi2}) and the
approximation~(\ref{approxsol}) give the definition
\begin{equation}\label{eq-defkappa} 
\begin{split}
  \kappa(t,x)& = - \frac{\int_{e_z>0 \text{ in } {\cal E}(0)}  e_z
    f^0(t,x,0,v_x,v_z(0,e_z)) \, de}{\int_{e_z<0 \text{ in } {\cal
        E}(0)} e_zG(e)\,  de} \\
& =\frac{\int_{v_z>0} \int v_z f^0(t,x,0,v_x,v_z)\, dv_xdv_z}{\int_{v_z>0} \int v_z M(v_x,v_z)\,  dv_xdv_z},
\end{split}
\end{equation}
where ${M}({v}_x,{v}_z)=G(v_x,e_z(0,v_z))=\exp (-({v}_x^2+{v}_z^2)/2)$. 

Now, note that~(\ref{orderzerointerfacevzpos}) has to be modified
according to the Knudsen layer correction to get
\begin{equation*}
\begin{split}
  {f}^0({t},{x},0,{v}_x,{v}_z)_{{v}_z<0} & =
  {\phi}^0({t},{x},0,{v}_x,{e}_z(0,v_z)) +
  {\psi}^0({t},{x},0,{v}_x,{e}_z(0,v_z)) \\
& = \chi({t},{x},0,{v}_x,{e}_z(0,v_z)). 
\end{split}
\end{equation*}
Consequently, the definition~(\ref{eq-defchi}) of $\chi$ and the
approximation~(\ref{approxsol}) give the following approximation of
the outgoing distribution
\begin{equation}\label{eq-approxoutg} 
{f}^0({t},{x},0,{v}_x,{v}_z)_{{v}_z<0} \approx  {\kappa}({t},{x}) {M}({v}_x,{v}_z), 
\end{equation}
which gives in dimensional variables the classical 
perfect accommodation boundary condition (\ref{perfaccomodationcond})
(sometimes called the diffuse reflexion boundary condition), provided
that the coefficient
$\kappa$ is such that the corresponding approximation of
the mass flux of $f^0$ at the boundary $z=0$ is zero. Indeed, the
definition~(\ref{eq-defkappa}) of $\kappa$ implies that this property holds.

\bigskip

%%%%% Maxwell %%%%
(3) Finally, we assume $\eta=O(1)$, which corresponds to $\tau_{ms}^*
\approx \tau_z^* \ll t_B^*$. The LKSL problem satisfied by $\phi^0$ is
\begin{eqnarray}
 {v}_z({z},{e}_z)\partial_{{z}}{\phi}^0 &=& {Q}[{\phi}^0],\; 
 \mbox{ for}\ {z}_+({e}_z)<{z}<{z}_-({e}_z),\  \label{orderzeromax-2}\\
{\phi}^0({t},{x},0,{v}_x,{e}_z)_{{e}_z>\sqrt{{W}_m}}&=&{f}^0({t},{x},0,{v}_x,{v}_z(0,{e}_z)),
\label{orderzerointerfacevznegmax-2}\\
{\phi}^0({t},{x},{z}_{-}({e}_z),{v}_x,{e}_z)&=& {\phi}^0({t},{x},{z}_{-}({e}_z),{v}_x,-{e}_z),
\label{orderzeroreflexionzm-2}\\
{\phi}^0({t},{x},{z}_{+}({e}_z),{v}_x,{e}_z)&=& {\phi}^0({t},{x},{z}_{+}({e}_z),{v}_x,-{e}_z), 
\mbox{ for}\; |{e}_z| < \sqrt{{W}_m}. \label{orderzeroreflexionzp}
\end{eqnarray}
We can claim that this linear kinetic surface layer (LKSL) problem has
a unique solution and that this solution has a zero mass flux through
the surface $z=0$ (see lemma 1 in the following): 
\begin{equation}
\int_{|{e}_z|>\sqrt{{W}_m}}\int {e}_z{\phi}^0({t},{x},0,{v}_x,{e}_z)\  d{v}_xd{e}_z=0. \label{orderzeromassflux}
\end{equation}

Now if we solve the LKSL problem (\ref{orderzeromax-2}-\ref{orderzeroreflexionzp}),
then $\phi^0(t,x,0,v_x,e_z)$, which is the value of the solution at ${z}=0$ for ${e}_z<0$, gives a boundary value
for~(\ref{orderzerointerfacevzpos}). This value linearily depends on
the inflow data: ${\phi}^0({t},{x},0,{v}_x,{e}_z(0,{v}_z))=
{ \cal A}\ ({f}^0({t},{x},0,{v}_x,.)_{|{v}_z>0})$, where $\cal A$ is called
the "albedo" operator $\cal A$. Consequently, the boundary
condition~(\ref{orderzerointerfacevzpos}) of~(\ref{eq-f0}) reads 
\begin{equation}
{f}^0({t},{x},0,{v}_x,{v}_z<0)=
{ \cal A}\ ({f}^0({t},{x},0,{v}_x,.)_{|{v}_z>0}). \label{exactBC}
\end{equation}
This relation can be interpreted as an exact boundary
condition. However, the operator ${\cal A}$ is implicitely defined: we must solve the LKSL problem
(\ref{orderzeromax-2}-\ref{orderzeroreflexionzp}) to get 
${\phi}^0({t},{x},0,{v}_x,{e}_z)_{{e}_z<0}$, which could be done approximately by a numerical computation.
Nevertheless, it is possible to get an approximation of the
operator ${\cal A}$ that explicitely gives
${\phi}^0({t},{x},0,{v}_x,{e}_z)_{{e}_z<0}$ as a function of
${f}^0({t},{x},0,{v}_x,{v}_z(0,{e}_z))_{|{e}_z>0}$: using again Lemma 1 we conclude that
\begin{eqnarray*}
{\phi}^0({t},{x},0,{v}_x,{e}_z)_{{e}_z<0}&\approx& {\phi}^{0,1}({t},{x},0,{v}_x,{e}_z)_{{e}_z<0}, \\
&\approx& 
(1-{a}({e}_z)){f}^0({t},x,0,{v}_x,-{v}_z(0,{e}_z) ) 
+ {a}({e}_z){\alpha} (t,x){G}({v}_x,{e}_z).
\end{eqnarray*}
From (\ref{orderzerointerfacevzpos}), we get
\begin{eqnarray*}
{f}^0({t},{x},0,{v}_x,{v}_z)_{{v}_z<0}&\approx& {\phi}^{0,1}({t},{x},0,{v}_x,{e}_z(0,v_z)) \\
&\approx& 
(1-{a}({v}_z)){f}^0({t},x,0,{v}_x,-{v}_z )  
+ {a}({v}_z)\beta(t,x)M({v}_x,{v}_z).
\end{eqnarray*}
Moreover, for the same reason as for the previous regime, the
approximation of mass flux of ${f}^0$ through the boundary ${z}=0$,
and hence the coefficient $\beta$ can be uniquely determined. Coming
back in dimensional variables, we get~(\ref{Maxcond}). From this
relation we can easily check that the associated scattering kernel
satisfies the properties of non-negativeness, normalization
and since $a(v)=a(-v)$, the property of reciprocity.\\

$\Box$

%%%%%%%%
%%% LEMMA 1
%%%%%%%%
\begin{lemma}
Let us consider the linear kinetic surface layer problem (LKSL) 
\begin{eqnarray}
  {v}_z({z},{e}_z)\partial_{{z}}{\phi}^0 &=& {Q}[{\phi}^0],\; 
 \mbox{ for}\ {z}_+({e}_z)<{z}<{z}_-({e}_z),\  \label{LKSL-1}\\
{\phi}^0(0,{v}_x,{e}_z)_{{e}_z>\sqrt{{W}_m}}&=&{f}^*({v}_x,{v}_z(0,{e}_z)),
\label{LKSL-2}\\
{\phi}^0({z}_{-}({e}_z),{v}_x,{e}_z)&=& {\phi}^0({z}_{-}({e}_z),{v}_x,-{e}_z),
\label{LKSL-3}\\
{\phi}^0({z}_{+}({e}_z),{v}_x,{e}_z)&=& \phi^0({z}_{+}({e}_z),v_x,-{e}_z), 
\mbox{ for}\; |{e}_z| < \sqrt{{W}_m}. \label{LKSL-4}
\end{eqnarray}
This problem has a unique solution 
%and thus defines a unique outgoing distribution function ${\phi}^0(0,v_x,e_z)$ for $e_z<-\sqrt{W_m}$ 
 and this solution has a zero mass flux through the surface $z=0$: 
\begin{equation}
\int_{|{e}_z|>\sqrt{{W}_m}}\int{e}_z{\phi}^0(0,{v}_x,{e}_z)\ d{v}_xd{e}_z=0.  \label{lkslflux}
\end{equation}
Moreover in a first approximation we have
\begin{eqnarray*}
{\phi}^0(0,{v}_x,{e}_z)_{|{e}_z<-\sqrt{{W}_m}}&\approx&
(1-{a}({e}_z)){f}^*({v}_x,-{v}_z(0,{e}_z) ) 
+ {a}({e}_z){\alpha}_1{G}({v}_x,{e}_z),
\end{eqnarray*}
where the coefficient $a$ is given by
\begin{equation}
a(e)=  1- \exp \left(-\frac{2{\tau}_z({e}_z)}{\bar{\tau}_{ms}({e})}\right) \label{coeffaccod}.
\end{equation}

\end{lemma}

{\em Proof :}\\
{\bf (i) Existence and uniqueness:} Existence and uniqueness of a
solution of the LKSL problem (\ref{LKSL-1}-\ref{LKSL-4}) can be proved
by using standard
techniques in linear transport problems. The reader can refer,  for
instance, to  \cite{Falk}. 
%Thus the outgoing distribution function is uniquely defined.
\\

\noindent
{\bf (ii) Mass flux at $z=0$ :}
Multiplying (\ref{LKSL-1}) by $|{e}_z|{\sigma}({z},{e}_z)$ and
using~(\ref{eq-sigmavzez}), we get
\[ {e}_z\partial_{{z}}{\phi}^0 =  {Q}[{\phi}^0]|{e}_z|{\sigma}({z},{e}_z) .\]
Now we integrate this relation with respect to ${z}$. It comes
\begin{eqnarray*}
 \int_{{z}_+({e}_z)}^{{z}_-({e}_z)} {e}_z\partial_{{z}}{\phi}^0\ d{z} =
 \int_{{z}_+({e}_z)}^{{z}_-({e}_z)}{Q}[{\phi}^0]|{e}_z|{\sigma}({z},{e}_z)\ d{z},
\end{eqnarray*}
or,
\begin{eqnarray*}
 {e}_z{\phi}^0({z}_-({e}_z),{v}_x,{e}_z)-
{e}_z{\phi}^0({z}_+({e}_z),{v}_x,{e}_z)= 
 \int_{{z}_+({e}_z)}^{{z}_-({e}_z)}{Q}[{\phi}^0]|{e}_z|{\sigma}({z},{e}_z)\ d{z},
\end{eqnarray*}
where ${z}_+({e}_z)=0$ for $|{e}_z|>\sqrt{{W}_m}$. Now integrating
with respect to ${v}_x$ and ${e}_z$, we find
\begin{equation}\label{eq-intLKSL} 
\begin{split}
\int \int {e}_z{\phi}^0({z}_-({e}_z),{v}_x,{e}_z)\ d{v}_xd{e}_z&-
\int \int{e}_z{\phi}^0({z}_+({e}_z),{v}_x,{e}_z)\ d{v}_xd{e}_z  \\
&= \int \int \int_{{z}_+({e}_z)}^{{z}_-({e}_z)}{Q}[{\phi}^0]|{e}_z|{\sigma}({z},{e}_z)\ d{z}d{v}_xd{e}_z. 
\end{split}
\end{equation}
But since ${e}_z{\phi}^0({z}_{\pm}({e}_z),{v}_x,{e}_z)\ d{v}_xd{e}_z$
is an odd function of ${e}_z$ (see~(\ref{LKSL-3}) and~(\ref{LKSL-4})),
the first term of the left-hand side of this relation vanishes and the
second one gives
\begin{eqnarray*}
\int \int{e}_z{\phi}^0({z}_+({e}_z),{v}_x,{e}_z)\ d{v}_xd{e}_z&=&
\int \int_{|{e}_z|<\sqrt{{W}_m}} {e}_z{\phi}^0({z}_+({e}_z),{v}_x,{e}_z)\ d{e}_zd{v}_x   \\
&+&\int\int_{|{e}_z|>\sqrt{{W}_m}}{e}_z{\phi}^0(0,{v}_x,{e}_z)\ d{e}_zd{v}_x, \\
&=& \int\int_{|{e}_z|>\sqrt{{W}_m}}{e}_z{\phi}^0(0,{v}_x,{e}_z)\ d{e}_zd{v}_x.
\end{eqnarray*}
Consequently,~(\ref{eq-intLKSL}) now reads
\begin{equation*}
-  \int\int_{|{e}_z|>\sqrt{{W}_m}}{e}_z{\phi}^0(0,{v}_x,{e}_z)\
d{e}_zd{v}_x = \int \int \int_{{z}_+({e}_z)}^{{z}_-({e}_z)}{Q}[{\phi}^0]|{e}_z|{\sigma}({z},{e}_z)\ d{z}d{v}_xd{e}_z.
\end{equation*}
Finally, inverting the integration with respect to ${z}$ and the integration
with respect to ${v}_x$ and ${e}_z$ in the right-hand side (see~(\ref{eq-fubini})), we get
\begin{eqnarray*}
-\int \int_{|{e}_z|>\sqrt{{W}_m}}{e}_z{\phi}^0(0,{v}_x,{e}_z)\ d{v}_xd{e}_z &=& 
 \int_{0}^{{L}} \int_{{\cal E}({z})}
 {Q}[{\phi}^0]|{e}_z|{\sigma}({z},{e}_z)\, de d{z}, \\
&=& 0,
\end{eqnarray*}
due to the mass conservation (see~(\ref{eq-cons})). \\

\noindent
{\bf (iii) Approximate solution of the LKSL problem} \\
First, we multiply~(\ref{LKSL-1}) by $\sigma(z,e_z)$, and we use the
decomposition of the collision operator into gain and loss terms to
rewrite~(\ref{LKSL-1}) as 
\begin{equation}\label{eq-eqphi0} 
    sign(e_z) \partial_{{z}}{\phi}^0 = \sigma(z,e_z){Q^+}[{\phi}^0]
    - \frac{{\sigma}({z},{e}_z)}{{\tau}_{ms}(z,{e}_z)}{\phi}^0, \mbox{ for}\ {z}_+({e}_z)<{z}<{z}_-({e}_z)
\end{equation}
Now we proceed by looking for
 an approximate solution of the boundary value problem (\ref{eq-eqphi0},\ref{LKSL-2}--\ref{LKSL-4}) in the form 
 ${\phi}^0={\phi}^{0,(1)}+{\phi}^{0,(2)}+ ....$, and we detail below
 how we construct the first approximation ${\phi}^{0,(1)}$.

Let us remark that if $f^*$ is a Maxwellian, then 
${\alpha}{G}({v}_x,{e}_z)$ is a solution of the LKSL problem for any
constant $\alpha$. Therefore, we propose to construct a first
approximation $\phi^{0,(1)}$ of $\phi^0$ as follows. We replace
$\phi^0$ in the gain term of~(\ref{eq-eqphi0}) by the Maxwellian
${\phi}^{0,(0)}= {\alpha_1}\ {G}({v}_x,{e}_z)$ (where the constant
$\alpha_1$ is undetermined for the moment) to get the following problem that defines 
%% Then, if $f^* \neq \alpha_1 G$,  we insert the approximate value 
%% ${\phi}^{0,(0)}$ in the operator $\cal T$, or which is equivalent 
$\phi^{0,(1)}$:
\begin{eqnarray}
  sign(e_z)\partial_{{z}}{\phi}^{0,(1)} &=& \frac{\sigma(z,e_z)}{\tau_{ms}(z,e_z)} \left({\alpha}_1\frac{{G}({e})}{{\tau}_{ms}(z,{e})}-
 {\phi}^{0,(1)} \right) , \quad  {z}_+({e}_z)<{z}<{z}_-({e}_z),\  \label{LKSL1-1}\\
{\phi}^{0,(1)}(0,{v}_x,{e}_z)_{{e}_z>\sqrt{{W}_m}}&=&{f}^*({v}_x,{v}_z(0,{e}_z)),
\label{LKSL1-2}\\
{\phi}^{0,(1)}({z}_{-}({e}_z),{v}_x,{e}_z)&=&
 {\phi}^{0,(1)}({z}_{-}({e}_z),{v}_x,-{e}_z),
\label{LKSL1-3}\\
{\phi}^{0,(1)}({z}_{+}({e}_z),{v}_x,{e}_z)&=& {\phi}^{0,(1)}({z}_{+}({e}_z),{v}_x,-{e}_z), 
\mbox{ for}\; |{e}_z| < \sqrt{{W}_m}. \label{LKSL1-4}
\end{eqnarray}
Then the solution $\phi^{0,(1)}$ can be explicitely constructed by integrating~(\ref{LKSL1-1}) along trajectories of free and trapped
molecules. This approach guarantees that the corresponding
approximated boundary condition~(\ref{Maxcond}) is exact is $f$ is a Maxwellian.\\

\noindent \underline{\em free molecules with $e_z>0$ ($e_z>\sqrt{W_m}$):} 

\smallskip

\noindent in that case, particles go from $z_+(e_z)=0$ to $z_-(e_z)$, and we can
integrate~(\ref{LKSL1-1}) between $0$ and some $z\in \lbrack
0,z_-(e_z) \rbrack$ to get
\begin{eqnarray}
{\phi}^{0,(1)}({z},{v}_x,{e}_z) &=& 
\exp \left(-\int_0^{{z}}\frac{{\sigma}({\zeta},{e}_z)}{{\tau}_{ms}(\zeta,{e})}d{\zeta}\right)
{f}^*({v}_x,{v}_z(0,{e}_z)) \nonumber \\
\nonumber &+& \exp \left(-\int_0^{{z}}\frac{{\sigma}({\zeta},{e}_z)}{{\tau}_{ms}(\zeta,{e})}d{\zeta}\right)
\int_0^{{z}} 
\exp \left(\int_0^{{\zeta}}\frac{{\sigma}({\eta},{e}_z)}{{\tau}_{ms}(\eta,{e})}d{\eta}\right)
\frac{{\sigma}({\zeta},{e}_z)}{{\tau}_{ms}(\zeta,{e})}d{\zeta}\ {\alpha}_1{G}.
\\ \label{phi01epos}
\end{eqnarray}
We write this relation at $z=z_-(e_z)$, we use the
definition of $\overline{\tau}_{ms}$ (see~(\ref{meantaums})), and then
exact computations of the integral of the exponential gives
\begin{equation}\label{eq-clfreepos} 
{\phi}^{0,(1)}({z}_{-}({e}_z),{v}_x,{e}_z)_{|{e}_z>0} 
= \exp \left(-\frac{{\tau}_z({e}_z)}{\overline{\tau}_{ms}({e})}\right)
{f}^*({v}_x,{v}_z(0,e_z)) +\left(1- \exp \left(-\frac{{\tau}_z({e}_z)}{\overline{\tau}_{ms}({e})}\right)  \right) {\alpha}_1{G}.
\end{equation}
% which is an even function of ${e}_z$ since $\overline{\tau}_{ms}({v}_x,{e}_z)=\overline{\tau}_{ms}({v}_x,-{e}_z)$.\\

\noindent \underline{\em free molecules with $e_z<0$ ($e_z<-\sqrt{W_m}$):} 

\smallskip

\noindent in that case, particles go from $z_-(e_z)$ to $z_+(e_z)=0$. First,
we use~(\ref{LKSL1-3}),~(\ref{eq-clfreepos}), the even parity
of $\bar{\tau}_{ms}$, $\tau_z$ and $G$ and the odd parity of
$v_z(0,e_z)$ with respect to $e_z$ to obtain the
distribution of outgoing particles at $z=z_-(e_z)$: 
\begin{equation*}
{\phi}^{0,(1)}({z}_{-}({e}_z),{v}_x,{e}_z)_{|{e}_z<0} 
= \exp \left(-\frac{{\tau}_z({e}_z)}{\overline{\tau}_{ms}({e})}\right)
{f}^*({v}_x,-{v}_z(0,e_z)) +\left(1- \exp \left(-\frac{{\tau}_z({e}_z)}{\overline{\tau}_{ms}({e})}\right)  \right) {\alpha}_1{G}.
\end{equation*}
Then we can
integrate~(\ref{LKSL1-1}) between $z_-(e_z)$ and some $z\in \lbrack
0,z_-(e_z) \rbrack$ to get 
\begin{eqnarray*}
&& {\phi}^{0,(1)}({z},{v}_x,{e}_z)=   \nonumber\\
&&\exp \left(-\int_{{z}}^{{z}_-}\frac{{\sigma}({\zeta},{e}_z)}{{\tau}_{ms}(\zeta,{e})}d{\zeta}\right)
\left( \exp \left(-\frac{{\tau}_z({e}_z)}{\overline{\tau}_{ms}({e})}\right){f}^*({v}_x,-{v}_z(0,{e}_z))  
+(1- \exp \left(-\frac{{\tau}_z({e}_z)}{\overline{\tau}_{ms}({e})}\right)){\alpha}_1{G}     \right) \nonumber \\
&&-\left(  \exp \left(-\int_{{z}}^{{z}_-}\frac{{\sigma}({\zeta},{e}_z)}{{\tau}_{ms}(\zeta,{e})}d{\zeta}\right)-1 \right) {\alpha}_1{G},
% \label{phi01enegF}
\end{eqnarray*}
and hence the distribution of outgoing particles at $z=0$:
\begin{eqnarray*}
{\phi}^{0,(1)}(0,{v}_x,{e}_z)_{|{e}_z<0}
% &=& \exp \left(-\frac{2{\tau}_z({e}_z)}{\overline{\tau}_{ms}({e})}\right)
% {f}^*({v}_x,-{v}_z(0,{e}_z))\\
% && +\exp \left(-\frac{{\tau}_z({e}_z)}{\overline{\tau}_{ms}({e})}\right)
% \left(1- \exp \left(-\frac{{\tau}_z({e}_z)}{\overline{\tau}_{ms}({e})}\right)  \right) {\alpha}_1{G} 
%  + \left(1- \exp \left(-\frac{{\tau}_z({e}_z)}{\overline{\tau}_{ms}({e})}\right)  \right) {\alpha}_1{G},\\
&=& \exp \left(-\frac{2{\tau}_z({e}_z)}{\overline{\tau}_{ms}({e})}\right) {f}^*({v}_x,-{v}_z(0,{e}_z))
+\left(1- \exp
  \left(-\frac{2{\tau}_z({e}_z)}{\overline{\tau}_{ms}({e})}\right)
\right) {\alpha}_1{G} . \label{eq-clmaxwphi01}
\end{eqnarray*}
This shows that we can construct the first approximation
$\phi^{0,(1)}$ for free particles, and that this approximation
satisfies a Maxwell boundary condition at $z=0$ with the accomodation
coefficient $a(e)=1- \exp
\left(-\frac{2{\tau}_z({e}_z)}{\overline{\tau}_{ms}({e})} \right)$, provided
  that the coefficient $\alpha_1$ can be defined such that the
  corresponding mass flux is zero. Indeed, using~(\ref{LKSL1-2})
  and~(\ref{eq-clmaxwphi01}), it is sufficient to set
\begin{equation*}
{\alpha}_1=\left(\int \int_{e_z>\sqrt{W_m}}e_z a(e)f^*({v}_x,-{v}_z(0,{e}_z))de_zdv_x\right)/ \left(  \int\int_{e_z>\sqrt{W_m}}e_za(e)G(v_x,e_z)de_zdv_x\right). 
\end{equation*}

Now, $\phi^{0,(1)}$ must also be constructed for trapped particles in
order to have a complete approximation of $\phi^{0}$. In the
following, we follow the same approach as that used for free particles.

\bigskip

\noindent \underline{\em trapped molecules with $e_z>0$ ($e_z<\sqrt{W_m}$):}

\smallskip
\noindent in that case, particles go from $z_+(e_z)$ to $z_-(e_z)$, and we can
integrate~(\ref{LKSL1-1}) between $z_+(e_z)$ and some $z\in \lbrack
z_+(e_z),z_-(e_z) \rbrack$ to get
\begin{eqnarray}
{\phi}^{0,(1)}({z},{v}_x,{e}_z) &=& 
\exp \left(-\int_{z_+(e_z)}^{{z}}\frac{{\sigma}({\zeta},{e}_z)}{{\tau}_{ms}(\zeta,{e})}d{\zeta}\right)
{\phi}^{0,(1)}({z}_+({e}_z),{v}_x,{e}_z) \nonumber \\
\nonumber &+& \exp \left(-\int_{z_+(e_z)}^{{z}}\frac{{\sigma} ({\zeta},{e}_z)}{{\tau}_{ms}(\zeta,{e})}d{\zeta}\right)
\int_{z_+(e_z)}^{{z}} 
\exp \left(\int_{z_+(e_z)}^{{\zeta}}\frac{{\sigma}({\eta},{e}_z)}{{\tau}_{ms}(\eta,{e})}d{\eta}\right)
\frac{{\sigma}({\zeta},{e}_z)}{{\tau}_{ms}(\zeta,{e})}d{\zeta}\ {\alpha}_1{G}, 
\\ \label{phi01epost}
\end{eqnarray}
and thus
\begin{equation}\label{eq-cltrappedpos} 
{\phi}^{0,(1)}({z}_{-}({e}_z),{v}_x,{e}_z)_{|{e}_z>0} 
= \exp \left(-\frac{{\tau}_z({e}_z)}{\overline{\tau}_{ms}({e})}\right)
{\phi}^{0,(1)}({z}_+({e}_z),{v}_x,{e}_z) +\left(1- \exp \left(-\frac{{\tau}_z({e}_z)}{\overline{\tau}_{ms}({e})}\right)  \right) {\alpha}_1{G}.
\end{equation}

\bigskip

\noindent \underline{\em trapped molecules with $e_z<0$ ($e_z>-\sqrt{W_m}$):}

\smallskip
\noindent in that case, particles go from $z_-(e_z)$ to $z_+(e_z)$. 
Consequently,
we use~(\ref{LKSL1-3}),~(\ref{eq-cltrappedpos}), the even parity
of $\bar{\tau}_{ms}$, $\tau_z$ and $G$ and the odd parity of
$v_z(0,e_z)$ with respect to $e_z$ to obtain the
distribution of outgoing particles at $z=z_-(e_z)$: 
\begin{equation*}
{\phi}^{0,(1)}({z}_{-}({e}_z),{v}_x,{e}_z)_{|{e}_z<0} 
= \exp \left(-\frac{{\tau}_z({e}_z)}{\overline{\tau}_{ms}({e})}\right)
{\phi}^{0,(1)}({z}_+({e}_z),{v}_x,-{e}_z) +\left(1- \exp \left(-\frac{{\tau}_z({e}_z)}{\overline{\tau}_{ms}({e})}\right)  \right) {\alpha}_1{G}.
\end{equation*}
Then we can
integrate~(\ref{LKSL1-1}) between $z_-(e_z)$ and some $z\in \lbrack
z_+(e_z),z_-(e_z) \rbrack$ to get
\begin{equation}\label{phi01enegT}
\begin{split}
  {\phi}^{0,(1)}({z},{v}_x,{e}_z) & =
\exp \left(-\int_{{z}}^{{z}_-}\frac{{\sigma}({\zeta},{e}_z)}{{\tau}_{ms}(\zeta ,{e})}d{\zeta}\right)
\left(  \exp
  \left(-\frac{{\tau}_z({e}_z)}{\overline{\tau}_{ms}({e})}\right){\phi}^{0,(1)}({z}_+({e}_z),{v}_x,-{e}_z)
\right. \\
& \left. \hspace{32ex} +(1- \exp \left(-\frac{{\tau}_z({e}_z)}{\overline{\tau}_{ms}({e})}\right)){\alpha}_1{G}     \right) \\
& \quad -\left(  \exp \left(-\int_{{z}}^{{z}_-}\frac{\ti{\sigma}({\zeta},{e}_z)}{{\tau}_{ms}(\zeta ,{e})}d{\zeta}\right)-1 \right) {\alpha}_1{G}, 
\end{split}
\end{equation}
and thus
\begin{equation*}
%\begin{split}
%& 
{\phi}^{0,(1)}({z}_+({e}_z),{v}_x,{e}_z)_{|{e}_z<0} %\\
% &\qquad = \exp \left(-\frac{2{\tau}_z({e}_z)}{\overline{\tau}_{ms}({e})}\right)
% {\phi}^{0,(1)}({z}_+({e}_z),{v}_x,-{e}_z)\\
% & \qquad \qquad  +\exp \left(-\frac{{\tau}_z({e}_z)}{\overline{\tau}_{ms}({e})}\right)\left(1- \exp \left(-\frac{{\tau}_z({e}_z)}{\overline{\tau}_{ms}({e})}\right)  \right) {\alpha}_1{G} 
%  + \left(1- \exp \left(-\frac{{\tau}_z({e}_z)}{\overline{\tau}_{ms}({e})}\right)  \right) {\alpha}_1{G},\\
%&\qquad 
= \exp \left(-\frac{2{\tau}_z({e}_z)}{\overline{\tau}_{ms}({e})}\right) 
{\phi}^{0,(1)}({z}_+({e}_z),{v}_x,-{e}_z)
+\left(1- \exp \left(-\frac{2{\tau}_z({e}_z)}{\overline{\tau}_{ms}({e})}\right)  \right) {\alpha}_1{G}.
%\end{split}
\end{equation*}
Then we can use~(\ref{LKSL1-4})
%But since
%${\phi}^{0,(1)}({z}_+({e}_z),{v}_x,-{e}_z)={\phi}^{0,(1)}({z}_+({e}_z),{v}_x,{e}_z)$
%we deduce
in the previous relation to deduce that
${\phi}^{0,(1)}({z}_+({e}_z),{v}_x,{e}_z)={\alpha}_1{G}$ for every
trapped particles. Finally, using~(\ref{phi01epost})
and~(\ref{phi01enegT}) we obtain
\begin{equation*}
{\phi}^{0,(1)}({z},{v}_x,-{e}_z)= {\alpha}_1{G},\; \mbox{ for}\ {z}_+({e}_z) \leq {z} \leq {z}_-({e}_z),
\end{equation*}
that is to say that trapped molecules are in equilibrium in
$\lbrack {z}_+({e}_z),{z}_-({e}_z) \rbrack$. The first approximation
$\phi^{0,(1)}$ of $\phi^0$ now is completely defined.

$\Box$

\section{Derivation of boundary conditions: wall with nanoscale roughness}
\label{sec:roughness}

\begin{figure}
 \centering
    \resizebox{!}{0.48\textheight} { \input{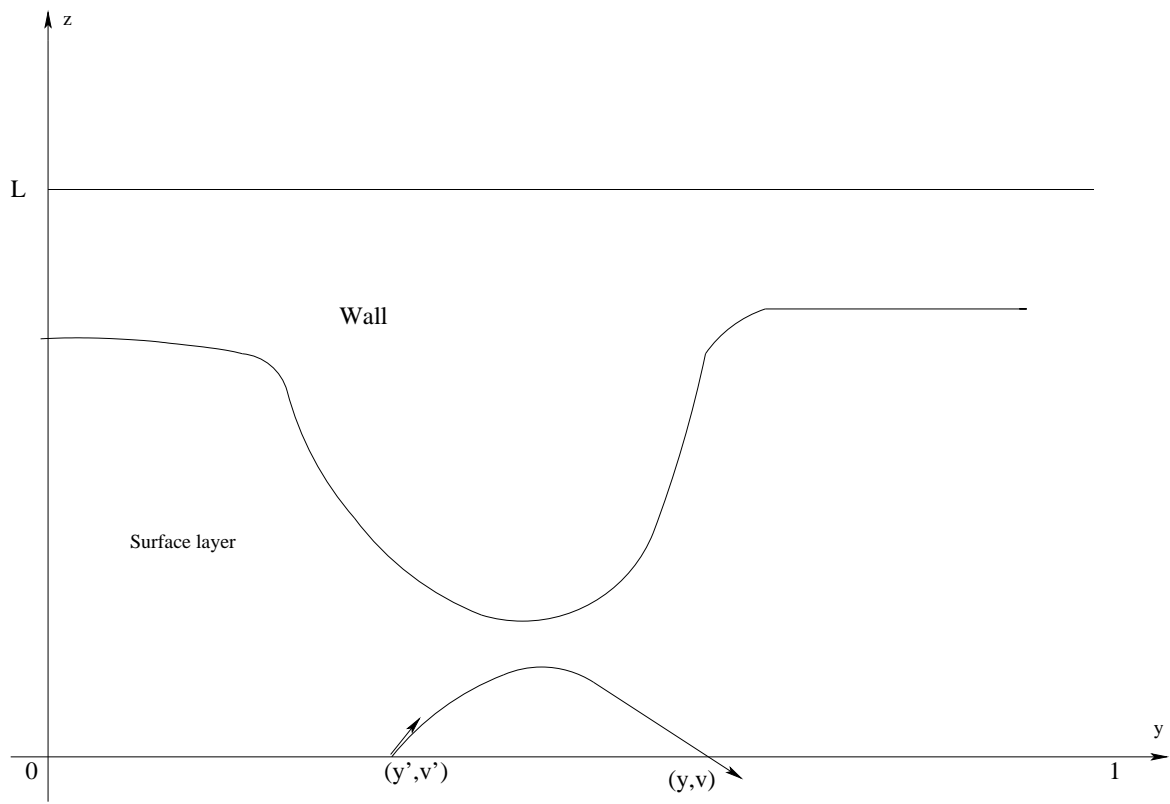}}
\caption{surface layer for a wall with nanoscale roughness and
 trajectory of a particle}
\label{fig:roughness}
\end{figure}

We assumed so far that the surface of the solid wall is flat and that the potential has the simplified form (\ref{potU+W}). 
Following the same approach, but with notations and algebra  a bit more tricky, we could obtain similar results for a more 
general attractive-repulsive surface potential ${\cal V}(x,z)$, corresponding to a smooth wall, i.e. such that ${\cal V}(x,z)=+\infty$ at $z=L$ .  
Moreover we can extend the approach to the case of a wall with nanoscale roughness (a wall on which there are a great number of minute asperities 
and which may induce multiple scattering
as indicated in \cite{Maxw}). More precisely,  let us consider the following configuration for the wall :  we assume that the surface layer is included 
in $[0, L]$ and that the potential ${\cal V}(x,z)$ is such that
\begin{equation}
{\cal V} (x,z)={\cal V}_{\#}(\frac{x}{L_*},z),  \label{periodicpot}
\end{equation}
where $L_*=\beta_*L$ and $\beta_*$ is a positive constant that
characterizes the roughness of the wall, and ${\cal V}_{\#}(y,z)$ is a periodic function of the nanoscopic variable 
$y$ with period $1$. This nanoscopic variable $y$ allows us to describe how a molecule impiging the surface layer at microscopic coordinate $x$ sees the 
nanoscopic roughness of the wall. 
  Moreover we assume that there exist $z=\zeta_{\infty} (y)$ a 1-periodic function with  $0\leq \zeta_{\infty} (y)<L$ and  $z=\zeta_0 (y)$ a 1-periodic function with 
$0\leq \zeta_0 (y) < \zeta_{\infty} (y)$ such that (see figure~\ref{fig:roughness})
\begin{eqnarray}
\lim_{z'<\zeta_{\infty} (y), z' \rightarrow \zeta_{\infty}(y)}{\cal V}_{\#}(y,z')&=&+\infty, \label{roughwall}\\
{\cal V}_{\#}(y,\zeta_0(y))&=&0. \label{roughwell}
\end{eqnarray}
 Finally, we assume that the potential is attractive-repulsive, i.e.
 \begin{equation}  
 \mbox{for }\; \zeta_0(y) <z<\zeta_{\infty}(y),\; \partial_z  {\cal V}_{\#}(y,z)>0, \;\;\;   \mbox{for }\; 0 <z<\zeta_0(y),\; \partial_z  {\cal V}_{\#}(y,z)<0,  \label{roughattractive-repulsive}
  \end{equation}
and that 
\begin{equation} {\cal V}_{\#}(y,z)= {\cal V}_m,\;\; \mbox{for}\; z\leq 0.  \label{roughfinite_range}
\end{equation}
The total energy of a molecule is 
\[E(x,z,v_x,v_z)=\frac{m}{2}|v|^2+{\cal V} (x,z), \]
and this total energy remains constant as long as the molecule does not collide with a phonon. 
Note that in this section, we do not use the change of velocity
variables $v\mapsto e(v,z)$. Indeed, since the potential is not
assumed to be separable into $U(x)+W(z)$ here, there is no obvious change of
variable that would simplify the equations.

With these assumptions, the flow of molecules is described by the following system of kinetic equations
\begin{equation}\label{eq-syst_rugosite} 
\begin{split}
& \partial_t f+v_x\partial_xf +v_z\partial_zf = 0, \; \; \;    z<0,   \\
& \partial_t f+v_x\partial_x f + v_z\partial_z f
-\frac{1}{m} \partial_x  {\cal V}(x,z)  \partial_{v_x}f
-\frac{1}{m} \partial_z  {\cal V}(x,z)\partial_{v_z}f= Q[f], \; \;
\;    0<z<L ,
\end{split}
\end{equation}
where the molecule-phonon collision term writes
\begin{equation}
Q[f]= \int  K(v,v') (M(v)f(v')-M(v')f(v))dv',
\end{equation}
and satisfies the properties recalled in proposition 1. Moreover the
distribution function $f$ is continuous through the interface $z=0$.

In the following, we compute the scattering kernel of asymptotic
boundary conditions corresponding to various regimes. However, we find
it more convenient to use the following form of the scattering kernel: 
\begin{equation*}
  k(v'\rightarrow v)=R(v'\rightarrow v)\frac{|v_z'|}{|v_z|},
\end{equation*}
where $R$ is the standard form (as used in~(\ref{eq-galbc})). 
%This new
%kernel is the probability density that a molecule reemitted with the
%velocity $v$ has striked the wall with a velocity between $v'$ and
%$v'+dv'$. 
With this new kernel, properties of normalization~(\ref{norm})
and reciprocity~(\ref{recip})
now read:
\begin{align}
& \int_{v_z<0 } k(v' \rightarrow v)v_z\ dv=-v_z', \label{norm_k} \\
& |v_z|k(v'\rightarrow v)M(v')=|v'_z|k(-v\rightarrow -v')M(v) . \label{recip_k}
\end{align}

%%%%%%%%%%%%%%%%%%%%%%%%%%%%%%%%%%%%%%%%%%%%%%%%%
%PROPOSITION 3
%%%%%%%%%%%%%%%%%%%%%%%%%%%%%%%%%%%%%%%%%%%%%%%%%
\begin{proposition}
Under the hypothesis (\ref{periodicpot}--\ref{roughfinite_range}), 
in the limit $\varepsilon =\frac{\tau_{ms}^*}{t_B^*}\rightarrow 0$,
the gas-surface interaction depends on the order of magnitude of
$\eta=\frac{\tau_{ms}^*}{\tau_{fl}^*}$ (where $\tau_{fl}^*$ is the characteristic
time of flight of a molecule through the surface layer), and can be described by the following boundary conditions at
$z=0$:
\begin{enumerate}
\item for $\eta=O(\frac{1}{\varepsilon})$, the boundary
  condition is the "specular" boundary condition which writes for a rough wall
\begin{equation}
f(t,x,0,v_x,v_z)_{| v_z<0} = \int_{ v'_z>0} k(v'\rightarrow v)f(t,x,0,v')\  dv', \label{roughspec}
\end{equation}
where the scattering kernel $k$, given by (\ref{roughspeckern}), 
is a probability density that 
is non-negative and satisfies the normalization and
reciprocity properties~(\ref{norm_k}--\ref{recip_k}).
\item for $\eta=O(1)$, the boundary condition writes, in a first
  approximation, as
\new{
\begin{equation}\label{roughmaxlike} 
f(t,x,0,v)= \int_{v'_z>0}k_1(v'\rightarrow v)f(t,x,0,v')\, dv' 
 + a^\#(v)\, \sigma(t,x)M(v)  ,
\end{equation}
where $k_1(v'\rightarrow v)$ can be viewed as a scattering kernel of non
thermalized molecules, and is defined by~(\ref{eq-k1}), $a^\#(v)$ is
the fraction of incident molecules that are re-emitted with the
velocity $v$ after a collision with a phonon (see~(\ref{eq-def_adiese})), and $\sigma$ is such
that the mass flux at $z=0$ is zero (defined in~(\ref{eq-sigma})). This boundary condition satisfies the properties of
non negativeness, normalization, and reciprocity.
}
\end{enumerate}
\end{proposition}

{\em Proof: } 
We denote by
$\phi=f_{|0<z<L}$, and we write $f$ 
and $\phi$ as functions of $(t,x,y=\frac{x}{L_*},z,v_x,v_z)$, periodic in $y$, with
period $1$. We use the same reference
quantities and nondimensional
variables as in the previous sections. With these new functions, the dimensionless form of
system~(\ref{eq-syst_rugosite}) is
 \begin{align}   
& \partial_t f+v_x\partial_xf
+\frac{\eta}{\beta_*\varepsilon}v_x\partial_yf+v_z\partial_zf = 0,
\quad \text{ for } z<0,    \label{dimlessboltz}  \\
& \partial_t \phi+v_x\partial_x \phi +\frac{\eta}{\beta_*\varepsilon}v_x\partial_y \phi+\frac{\eta}{\varepsilon}v_z\partial_z \phi
-\frac{1}{2\beta_*}\frac{\eta}{\varepsilon} \partial_y  {\cal V}_{\#}(y,z)  \partial_{v_x}\phi
-\frac{1}{2}\frac{\eta}{\varepsilon}  \partial_z  {\cal V}_{\#}(y,z)\partial_{v_z}\phi= \frac{1}{\varepsilon}Q[\phi],
\end{align}
for $0<z<1$, with interface conditions
\begin{eqnarray}
\phi(t,x,y,0,v_x,v_z)_{v_z>0}&=&f(t,x,y,0,v_x,v_z),  \label{LKSL-r-2-adim_pre}\\
  f(t,x,y,0,v_x,v_z)_{|v_z<0}  &=&     \phi(t,x,y,0,v_x,v_z) \label{boltzmann-r-2-adim_pre}
\end{eqnarray}
for every $y \in  [0,1]$. 

We define the average of $f$ over a period: 
\begin{equation*}
  F(t,x,z,v_x,v_z)=\int_0^1 f(t,x,y,z,v_x,v_z)\, dy.
\end{equation*}
Integrating (\ref{dimlessboltz}) with
respect to $y$ and taking into account the 1-periodicity , we obtain
\begin{equation}
\partial_t F +v_x\partial_x F +v_z\partial_z F  = 0,   \label{meanboltz}
\end{equation}
for $z<0$, and the boundary condition~(\ref{boltzmann-r-2-adim_pre})
leads to
\begin{equation} \label{boltzmann-r-2-adim}
%\phi(t,x,y,0,v_x,v_z)_{v_z>0}&=&f(t,x,y,0,v_x,v_z),\quad \forall  y \in  [0,1],  \label{LKSL-r-2-adim}\\
 F(t,x,0,v_x,v_z)_{|v_z<0}  =    \int_0^1  \phi(t,x,y,0,v_x,v_z)dy.
\end{equation}
Now, we use an expansion of $f$, $F$, and $\phi$ in terms of powers of
$\varepsilon$, and we identify the terms of same order of
magnitude. \\

(1) We consider the case $\eta=O(\frac{1}{\varepsilon})$, which
implies that $\tau_{fl}^*\ll \tau_{ms}^* \ll t_B^*$ (i.e. a weak
molecule-phonon interaction). 
%  \begin{eqnarray}
% \partial_t f+v_x\partial_xf +\frac{v_x}{\varepsilon^2}\partial_yf+v_z\partial_zf &=& 0, \; \; \;    z<0,    \label{dimlessboltz} \\
% \partial_t \phi+v_x\partial_x \phi +\frac{v_x}{\varepsilon^2}\partial_y \phi+\frac{v_z}{\varepsilon^2}\partial_z \phi
% -\frac{1}{\varepsilon^2} \partial_y  {\cal V}_{\#}(y,z)  \partial_{v_x}\phi
% -\frac{1}{\varepsilon^2} \partial_z  {\cal V}_{\#}(y,z)\partial_{v_z}\phi&=& Q[\phi], \; \; \;    0<z<L . \nonumber
% \end{eqnarray}
We find at zeroth order $\partial_yf^0= 0$ for $z<0$, which means
that $f^0$ does not depend on $y$, and
hence $F^0(t,x,z,v_x,v_z)=f^0(t,x,z,v_x,v_z)$. Consequently,
equations~(\ref{meanboltz}) and~(\ref{boltzmann-r-2-adim}) give 
\begin{eqnarray}
\partial_tF^0+v_x\partial_x F^0 +v_z\partial_zF^0  &=& 0, \; \; \;    z<0,  \label{meanboltz-2} \\
F^0(t,x,0,v_x,v_z)_{v_z<0}&=& \int_0^1  \phi^0(t,x,y,0,v_x,v_z)dy. \label{boltzmann-r-2}
\end{eqnarray}
However, $\phi^0$ stil depends on $y$ and we get 
\begin{equation}   \label{LKSL-r-1}
\frac{1}{\beta_*}v_x\partial_y \phi^0+ v_z\partial_z \phi^0
- \frac{1}{2\beta_*}\partial_y  {\cal   V}_{\#}(y,z)  \partial_{v_x}\phi^0
-\frac{1}{2} \partial_z  {\cal V}_{\#}(y,z)\partial_{v_z}\phi^0=0,
\end{equation}
with a boundary condition coming from~(\ref{LKSL-r-2-adim_pre}) which is 
\begin{equation}  \label{LKSL-r-2}
\phi^0(t,x,y,0,v_x,v_z)_{v_z>0}=F^0(t,x,0,v_x,v_z),\; \forall  y \in  [0,1].
\end{equation}
Note that the zeroth-order
system~(\ref{meanboltz-2}--\ref{LKSL-r-2}) in $(F^0,\phi^0)$ is
closed, contrary to the original system~(\ref{dimlessboltz}),
(\ref{LKSL-r-2-adim_pre}), (\ref{meanboltz}),
(\ref{boltzmann-r-2-adim_pre}) in $(F,\phi)$.

Relation (\ref{LKSL-r-2}) means that the molecules impinging the
surface layer with velocity $v=(v_x,v_z)$ see the roughness of the
wall from any nanoscopic variable $y$ with the same
probabiliy. Relation (\ref{boltzmann-r-2}) means that the number of
molecules going out of the surface layer at microscopic point $x$ with
velocity $v$ is the sum over $y$ of molecules going out with velocity
$v$ at the nanoscopic points $y,\; y \in [0,1]$. 

The characteristic curves of the LKSL problem (\ref{LKSL-r-1}),
defined by $\dot{y}(t)=v_x(t)/\beta_*$, $\dot{z}(t)=v_z(t)$,
$\dot{v}_x(t)=-\partial_y{\cal V}_{\#}(y(t),z(t))/2\beta_*$, and
$\dot{v}_z(t)=-\partial_z{\cal V}_{\#}(y(t),z(t))/2$, are the
trajectories of the molecules in the surface potential field. We
denote by $(y,v)=(y(y',v'),v(y',v'))=\Lambda (y',v')$ the mapping
that gives the position and the velocity $(y,v)$ of a molecule leaving the
surface layer (i. e. with $v_z<0$ at $z=0$) as a function of its position and
velocity $(y',v')$ when entering the surface layer (i.e. with
$v'_z>0$ at $z=0$), see figure~\ref{fig:roughness}.  Note that due to the time
reversibility of these trajectories, we have the important property
\new{\begin{equation}  \label{eq-charact}
(y,v)=\Lambda(y',v') \quad \Leftrightarrow \quad (y',-v')=\Lambda(y,-v),
\end{equation}
and hence $v'=-\Lambda_2(y,-v)$ for every
$(y,v,y',v')$ related by a characteristic curve. Another important
property is that the Jacobian of the transformation
$(y,v)=\Lambda(y',v')$ can be computed so that we have: 
\begin{equation}  \label{eq-jac}
|v'_z|dy'dv'=|v_z|dydv,
\end{equation}
see a proof in appendix~\ref{app:chvar}. The last property is 
that the total energy is
conserved along the characteristic and the potential energy has the
same value ${\cal V}_m$ at the head $(y,0)$ and the foot $(y',0)$ of
this characteristic, which yields
\begin{equation}  \label{eq-charnorme}
|\Lambda_2(y',v')|=|\Lambda_2(y,-v)|=|v|=|v'|.
\end{equation}
These relations are
essential to derive a collision kernel for
problem~(\ref{meanboltz-2}-\ref{boltzmann-r-2}) and to prove some of
its properties.}

Let $y$ in $[0,1]$ and $v$ such that $v_z<0$. Then
using the fact that $\phi^0$ is constant along the characteristics, we
get  $\phi^0(t,x,y,0,v)=\phi^0(t,x,y',0,v')$, where $(y',v')$ are such
that $(y,v)=\Lambda(y',v')$. Then using~(\ref{LKSL-r-2}) and the previous relation
$v'=-\Lambda_2(y,-v)$, we get
\begin{equation*}
  \phi^0(t,x,y,0,v)=F^0(t,x,0,-\Lambda_2(y,-v)).
\end{equation*}
Finally, we inject this relation into~(\ref{boltzmann-r-2}) to get
\begin{equation*}
F^0(t,x,0,v_x,v_z)_{v_z<0}= \int_0^1  F^0(t,x,0,-\Lambda_2(y,-v))\, dy,
\end{equation*}
which can be rewritten
\begin{equation}\label{roughspec2}
  F^0(t,x,0,v_x,v_z)_{v_z<0}= \int_{v'_z>0} k(v'\rightarrow v)  F^0(t,x,0,v')\, dv',
\end{equation}
where the collision kernel $k$ is defined by
\begin{equation}  \label{roughspeckern}
k(v'\rightarrow v)= \int_0^1  \delta(v'+\Lambda_2(y,-v))\, dy.
\end{equation}
This kernel is obviously non-negative, and it satisfies
\begin{equation*}
\int_{v'_z>0 }k(v'\rightarrow v) \ dv'=1,
\end{equation*}
and hence is a probability density.  
Indeed, note that a direct
integration of~(\ref{roughspeckern}) with respect to $v'$ and the use
of variables $v,y$ give this result. 

\new{The normalization property~(\ref{norm_k}) is obtained as follows:
  first, we use~(\ref{roughspeckern}) to get
\begin{equation*}
\int_{v_z<0}k(v'\rightarrow v)v_z \,dv
=-\int_{v_z<0}\int_0^1\delta(v'+\Lambda_2(y,-v))|v_z|\, dy dv.
\end{equation*}
Then, we use the change variables
  $(y,v)=\Lambda(y',w')$ and its properties~(\ref{eq-charact})
  and~(\ref{eq-jac}) to get
\begin{equation*}
\begin{split}
\int_{v_z<0}k(v'\rightarrow v)v_z \,
& = -\int_{w'_z>0}\int_0^1 \delta(v'-w')|w'_z|\, dy'dw'\\
& = -\int_0^1|v'_z|\, dy'  = -v'_z.
\end{split}
\end{equation*}
}

\new{Finally, the reciprocity property~(\ref{recip_k}) is obtained as
follows. First, we consider a given velocity $v$ (with $v_z<0$) and a
test function $\theta$, and we use~(\ref{norm_k}) to get
\begin{equation} \label{eq-reck1}
\begin{split}
\int_{v'_z>0}|v_z|k(v'\rightarrow v)M(v')\theta(v')\, dv' & 
=\int_{v'_z>0}\int_0^1|v_z|\delta(v'+\Lambda_2(y,-v))M(v')\theta(v')\,dydv' \\
& = \int_0^1|v_z|M(\Lambda_2(y,-v))\theta(-\Lambda_2(y,-v))\, dy\\
& = \int_0^1\theta(-\Lambda_2(y,-v))\, dy \, |v_z|M(v)\\
\end{split}
\end{equation}
from~(\ref{eq-charnorme}). Moreover,~(\ref{norm_k}) also gives
\begin{equation*}
\begin{split}
\int_{v'_z>0}|v'_z|k(-v\rightarrow -v')M(v)\theta(v')\, dv' & 
=\int_{v'_z>0}\int_0^1|v'_z|\delta(-v+\Lambda_2(y,v'))M(v)\theta(v')\,dydv'.
\end{split}
\end{equation*}
Then, we write $y'$ instead of $y$, and we use the change of variables
$(y,w)=\Lambda(y',v')$ and its properties~(\ref{eq-charact})
and~(\ref{eq-jac}) to get
\begin{equation} \label{eq-reck2}
\begin{split}
\int_{v'_z>0}|v'_z|k(-v\rightarrow -v')M(v)\theta(v')\, dv' & 
=\int_{w_z<0}\int_0^1|w_z|\delta(-v+w)M(v)\theta(-\Lambda_2(y,-w))\,dydw
\\
& =\int_0^1\theta(-\Lambda_2(y,-v))\,dy \, |v_z|M(v).
\end{split}
\end{equation}
Then, we compare~(\ref{eq-reck1}) and~(\ref{eq-reck2}) to find that
the two left-hand sides are equal for every test function $\theta$. The
reciprocity property $|v_z|k(v'\rightarrow
v)M(v')=|v'_z|k(-v\rightarrow -v')M(v)$ follows.
}

\bigskip

(2) Now we consider the case $\eta=O(1)$, which means that
$\tau_{fl}^*$ (the characteristic time of flight of a molecule across the surface layer) is 
comparable with $\tau_{ms}^*$ (the characteristic time of
molecule-phonon relaxation). \new{The zeroth order terms of the
  expansion are still denoted by $F^0$ and $\phi^0$, where $F^0$
  satisfies the same equation: 
\begin{align}
& \partial_tF^0+v_x\partial_x F^0 +v_z\partial_zF^0  = 0, \; \; \;    z<0,  \label{meanboltz-2-bis} \\
& F^0(t,x,0,v_x,v_z)_{v_z<0}= \int_0^1  \phi^0(t,x,y,0,v_x,v_z)dy, \label{boltzmann-r-2-bis}
\end{align}
and $\phi^0$ now
  is the periodic solution of 
\begin{align}
& \frac{1}{\beta_*}v_x\partial_y \phi^0+ v_z\partial_z \phi^0-
\frac{1}{2\beta_*}\partial_y  {\cal
  V}_{\#}(y,z)  \partial_{v_x}\phi^0-\frac{1}{2}\partial_z  {\cal
  V}_{\#}(y,z)\partial_{v_z}\phi^0=Q[\phi^0], 
  \label{LKSL-r-1-bis}\\
& \phi^0(t,x,y,0,v_x,v_z)_{v_z>0}=F^0(t,x,0,v_x,v_z),\; \forall  y \in
[0,1].   \label{LKSL-r-2-bis}
\end{align}
}
As in section 3.1, the right-hand-side of this equation is approximated
by $Q_+[\alpha (t,x){\cal M}]-\frac{\phi^0}{\tau_{ms}}$, where
$\tau_{ms}(v)=( \int K(v,v') M(v')\, dv')^{-1}$ is the
molecule-phonon relaxation time, and where ${\cal
  M}(y,z,v_x,v_z)=\exp (-|v|^2/2-{\cal V}_{\#}(y,z))$ which is
constant along the characteristics\new{, and $\alpha$ is a free parameter that
will be determined later}. To
integrate~(\ref{LKSL-r-1-bis}), it is useful to define the mean
molecule-phonon relaxation time $\bar{\tau}_{ms}(y',v')$ along the characteristic curve passing
by $(y',0, v')$ by
$\bar{\tau}_{ms}(y',v')=(\frac{1}{\tau_{fl}(y',v')}\int_{0}^{\tau_{fl}(y',v')}\frac{1}{\tau_{ms}(v(s))}\,
ds)^{-1}$. Then the solutions of
(\ref{LKSL-r-1-bis}) with boundary condition~(\ref{LKSL-r-2-bis}) satisfy
\begin{equation}\label{eq-sol_LKSL_approx} 
\begin{split}
\phi^0(t,x,y,0,v)=&\exp (-r(y',v'))F^0(t,x,0,v')\\
&+ (1-\exp (-r(y',v')))\sigma(t,x)M(v),
\end{split}
\end{equation}
where $r(y',v')=  \tau_{fl}(y',v')/\bar{\tau}_{ms}(y',v')$, $
\tau_{fl}(y',v')$ is the free time of flight of a molecule across the
surface layer in which it enters at $(y',z=0,v'))$, and \new{$\sigma(t,x)=\alpha(t,x)\exp(-{\cal
  V}_m)$ is still to be determined. 
First, note that $r(y',v')=r(y,-v)$: indeed it is defined as the
ratio of the free time of flight of a molecule along the trajectory
that starts at $(y',v')$ and ends at $(y,v)$ and the mean relaxation
time along this trajectory. Since this trajectory is the same as the
one that starts at $(y,-v)$ and ends at $(y',-v')$
(see~(\ref{eq-charact})), these two times are the same at $(y',v')$ and
$(y,-v)$. Then~(\ref{eq-sol_LKSL_approx}) can be rewritten as
\begin{equation}\label{eq-sol_LKSL_approx_ter} 
\begin{split}
\phi^0(t,x,y,0,v)=&\exp (-r(y,-v))F^0(t,x,0,-\Lambda_2(y,-v))\\
&+ (1-\exp (-r(y,-v)))\, \sigma(t,x)M(v),
\end{split}
\end{equation}
where $v'$ has been replaced by $-\Lambda_2(y,-v)$ due
to~(\ref{eq-charact}).
}

\new{Now, we use~(\ref{boltzmann-r-2-bis}), and the outgoing
  distribution $F^0(t,x,0,v)$ is found to be
\begin{equation}\label{eq-routgoing} 
\begin{split}
F^0(t,x,0,v)=
& \int_0^1\left(\exp (-r(y,-v))F^0(t,x,0,-\Lambda_2(y,-v)) 
+ (1-\exp (-r(y,-v)))\, \sigma(t,x)M(v)  \right)\, dy \\
 = & \int_0^1\int_{v'_z>0}(\exp(-r(y,-v))F^0(t,x,0,v')\delta(v'+\Lambda_2(y,-v)) \, dv'dy \\
& + \int_0^1 (1-\exp (-r(y,-v)))\, \sigma(t,x)M(v)  )\, dy \\
& = \int_{v'_z>0}k_1(v'\rightarrow v)F^0(t,x,0,v')\, dv' 
 +  \left(1-\int_0^1\exp (-r(y,-v))\, dy \right)\, \sigma(t,x)M(v)  ,
\end{split}
\end{equation}
where $k_1(v'\rightarrow v)$ can be viewed as a scattering kernel of non
thermalized molecules, and is defined by
\begin{equation}   \label{eq-k1}
k_1(v'\rightarrow v)=\int_0^1\exp (-r(y,-v))\delta(v'+\Lambda_2(y,-v)) \, dy.
\end{equation}
Moreover, the coefficient $a^\#(v)$ of~(\ref{roughmaxlike}) is found
to be 
\begin{equation}\label{eq-def_adiese} 
a^\#(v)=  1-\int_0^1\exp (-r(y,-v))\, dy .
\end{equation}
The computation of $\sigma$ and 
the reciprocity of this boundary condition are proved in appendix~\ref{app:recip_maxwell}.}

\section{Comments and concluding remarks}
\label{sec:conclusion}

1- In this approach, the boundary for the Boltzmann equation is considered to be located at $z=0$ which is the outer limit of the surface layer. 
The surface layer is considered as belonging to the solid phase. This is a two-phase description in opposition to the nanoscale models which are 
three-phase models (gas, surface layer, solid).\\

\noindent
2- The boundary condition (\ref{Maxcond}) is a Maxwell-like condition
but the "accommodation coefficient" $a=a(v)$ depends on the velocity.
More precisely the coefficient $a(v)$ must be interpreted as the
fraction of diffusively evaporated molecules.  A Maxwell-like
condition with a coefficient depending on the velocity has been
previously given in \cite{Borman}. Nevertheless the authors propose a
different expression : $\hat{a}(v)=
\frac{1}{1+(\tau_{ms}/(2\tau_z))}$. Let us remark that $\hat{a}(v)$
can be interpreted as a Pade approximant of $a(v)$ given in
(\ref{accomodation}), which can be explained since the boundary
condition is derived in \cite{Borman} from a nanoscale kinetic model
obtained by averaging (\ref{boltzmannmax}-\ref{reflexion2}) over the
surface layer. \\
 
  \noindent
  3- It is classical in the litterature (see for instance \cite{C}) to introduce the so-called accommodation
  coefficients $\alpha (\varphi )$ to describe the the interaction of a gas with a surface
\begin{equation*}
  \alpha (\varphi ) = \frac{\int_{v_z>0}\int |v_z| \varphi (v) \phi(v) dv_xdv_z-\int_{v_z<0}\int |v_z| \varphi (v) \phi(v) dv_xdv_z}{\int_{v_z>0}\int |v_z| \varphi (v) \phi(v) dv_xdv_z-J_0\int_{v_z<0}\int |v_z| \varphi (v) M(v) dv_xdv_z},
\end{equation*}
where $J_0=\int_{v_z>0}\int v_za(v_z){\phi}(0,{v}_x,{v}_z)dv_xdv_z/ \int_{v_z>0}\int v_za(v_z)M({v}_x,{v}_z)dv_xdv_z$,
and $\varphi (v)= v_x$ or $ v_z$, or $ |v|^2/2$ (accommodation coefficient for tangential or normal momentum or for
energy).  A drawback of the Maxwell's boundary condition noted in \cite{C} is that those various accommodation
coefficients are equal, and equal to the factor $a$ ( which explains why this coefficient is often called the
accommodation coefficient), which is not realistic since it is well-known that momentum and energy accommodate
differently in physical interactions. In contrast,
the boundary condition (\ref{Maxcond}) derived in the present paper gives different accommodation coefficients for energy and momentum. \\

\noindent
4- We notice that the boundary conditions obtained by this approach do not contain any free parameter to be
adjusted. All the information comes from the smaller scale (nanoscale). In particular the coefficient $a$ in the
Maxwell-like condition is given provided the interaction potential is known (and thus, $\tau_z$ and $\tau_{ms}$).  It is
interesting to look at the influence of the velocity on the fraction of diffusively evaporated molecules $a(v)$. Since
we assumed that the scattering kernel of the molecule-phonon collision term is bounded below and above (\ref{Kbound}),
then so is $\tau_{ms}$. Thus the behavior of $a$ for large $|v|$ depends essentially on $\tau_z$. Since
$\lim_{|v_z|\rightarrow +\infty}\tau_z(v_z)=0$, it appears that the fraction of diffusively evaporated molecules tends
to decrease for high velocities. Finally we remark also that a perfect accommodation boundary condition can be obtained
even if the interaction potential is purely repulsive. This is in contradiction with the
idea that a diffusive departure of molecules from a surface is due to desorption of trapped molecules (see \cite{Borman}).\\

\noindent
5- We considered in section~\ref{sec:roughness} rough walls with a periodicity assumption.  This assumption allows to
take into account the roughness of a surface in a simple way. Such a technique is commonly used for molecular dynamics
simulations in gas-surface interaction or in related applications such as porous media.  It can be relevant for instance
when the solid is a crystal or a composite material. Of course, for a wall with nanoscale roughness, even
the "specular" reflexion condition depends on the description of the surface potential.
The smooth wall can be seen as a particular case of a rough
wall by taking $k(v' \rightarrow v)=\delta_{v-(v'_x,-v'_z)}$.
Of course, taking advantage of such models for practical numerical
simulations requires accurate experiments to characterize the various
parameters for a given material. For computational purposes, an
approximation of the scattering kernel $k(v_p\rightarrow v_q)$ (for
$(v_p,v_q)$ in a discrete velocity grid) can be obtained by numerical
solutions of the characteristic curves in a unit cell of the surface layer.
\\

\noindent
6- In any boundary condition, the population of trapped molecules is not taken into account. This is justified when we
consider a bulk flow in a domain whose size is much greater than the thickness $L$ of the boundary layer. But when the
size of the domain becomes smaller (for instance in a channel with diameter comparable with $L$), then this might not be
correct. In such a configuration the number of molecules trapped in the surface layer cannot be neglected. Indeed if
we assume that the flow is stationary, then the distribution function writes
\[ \frac{m\ n_0 }{2k \pi T}e^{-{\cal V}(x,z)/kT}e^{-m(v_x^2+v_z^2)/2kT}.\]
Thus, the ratio of the number density of gas molecules at the outer boundary of the surface layer (and in the channel, i.e at $z\leq 0$)
over the number density of gas molecules at the bottom of the well potential (i.e at $z=z_*$)   
 is equal to
\[\frac{n(x,0)}{n(x,z_*)}= e^{-{\cal V}_m/kT},  \]
so that the number density of gas molecules inside the surface layer is much larger that the number density of gas
molecules in the channel when $kT\ll{\cal V}_m$ (see for instance \cite{KBA} for numerical results by means of molecular
dynamics simulations).  Thus in the vicinity of the wall we have to take into account the molecules inside the surface
layer, for instance to estimate the
mass flux parallel to the wall.\\

To conclude this paper we recall how the gas-surface interaction is described by the proposed kinetic approach (for a
smooth wall) at different scales and for various regimes.
\begin{itemize}
\item At the smaller scale
(the nanoscale, i.e. on a domain $[0,x^*] \times [0,z^*]$ with $x^*$ and $z^*$ $\approx$ 1 nanometer), the gas interaction is 
described by the two-dimensionnal kinetic model for the flow inside the surface layer 
(\ref{boltzmannmax}--\ref{reflexion2}),
suggested in \cite{Borman} and \cite{ACD} (coupled with the Boltzmann equation for the bulk flow). Then the gas-solid interaction 
at larger scales is derived from this model by formal systematic asymptotic analysis with various convenient scalings.

\item If we consider a gas flow in a domain $[0,x^{**}] \times
  [0,z^{**}]$ with $x^{**}\approx z^{**}\gg$ 1 nanometer, but where
  $x^{**}$ is the characteristic length of evolution of the flow in
  the x-direction inside the surface layer, then the gas-surface
  interaction can be described by the Boltzmann equation coupled with
  a one-dimensional kinetic or diffusion model describing the flow
  inside the surface layer of adsorbed molecules (mobile adsorption)
  which can be interpreted as non local boundary conditions for the
  Boltzmann equation in the bulk flow (see~\cite{ACD}).

\item At a larger scale, we consider a gas flow on a domain
  $[0,x^{***}] \times [0,z^{***}]$ where $x^{***}\approx z^{***}$ is
  the characteristic length of evolution of the Boltzmann equation in
  the bulk flow.  Then the gas-surface interaction
  can be described by the Boltzmann equation coupled with a local
  boundary condition that depends on the ratio $\tau_{ms}^*/\tau_z^*$:
\begin{itemize}
\item If $\tau_{ms}^*$, the characteristic time of relaxation of the
  molecules by the phonons, and $\tau_z^*$, the characteristic time
  for a molecule to cross the surface layer are comparable, then this
  boundary condition is implicitly given through the solution of a
  one-dimensional boundary value problem for a linear transport
  equation.  This boundary condition can be approximated by the
  numerical solution of the boundary value problem but it can
  also %when $\tau_{ms}^*/\tau_z^*$ is small this boundary condition can
  be approximated at first order by a Maxwell-like condition with a
  factor (the fraction of diffusively evaporated molecules) that
  depends on the velocity of the molecules. This fraction also depends
  on the temperature of the wall (through $M$, $\bar{\tau}_{ms}$, and
  $\bar{r}$) and of the morphology of the surface (through $\Lambda$).
\item If $\tau_z^ *\ll\tau_{ms}^*$, then the local boundary condition
  obtained is the well-known specular reflexion. 
\item If when $\tau_{ms}^ *\ll \tau_z^*$, then the local boundary
  condition obtained is the classical perfect accommodation boundary condition.
\end{itemize}

\end{itemize}

%%%%%%%%%%%%%%%%%%%% figures

\newpage

\clearpage

\newpage

\clearpage

\newpage

\clearpage
\clearpage

\newpage
\appendix
\section{Computation of the Jacobian of the transformation
  $(y,v)=\Lambda(y',v')$}
\label{app:chvar}
Let $\varphi_{in} (y',v')$ be a function defined for $y'\in [0,1]$ and
$v'$ such that $v'_z>0$. Let $\varphi_{out}(y,v)$ be the outgoing
value (at $z=0$, for $v_z<0$) of the solution $\varphi(y,z,v)$
of~(\ref{LKSL-r-1}), with $\varphi_{in} (y',v')_{|v'_z>0}$ as an
inflow boundary data (at $z=0$), \new{that is to say, $\varphi$ is a
  $y$-periodic function solution of: 
\begin{equation}   \label{LKSL-r-1-app}
\begin{split}
& \frac{1}{\beta_*}v_x\partial_y \varphi+ v_z\partial_z \varphi
- \frac{1}{2\beta_*}\partial_y  {\cal   V}_{\#}(y,z)  \partial_{v_x}\varphi
-\frac{1}{2} \partial_z  {\cal V}_{\#}(y,z)\partial_{v_z}\varphi=0, \\
& \varphi(y',0,v')_{v'_z>0}=\varphi_{in}(y',v') ,\; \forall  y' \in  [0,1].
\end{split}
\end{equation}
}
Since this solution is constant along
the characteristics, we have:
\begin{equation}
\varphi_{out} (y,v)=\varphi_{in} (y',v'), \label{carac}
\end{equation}
where $(y,v)=\Lambda (y',v')$ has been introduced
before~(\ref{eq-charact}). 

Denote by $J$ the Jacobian of the change of variables $(y,v)=\Lambda
(y',v')$, then, using~(\ref{carac}), we can write the average outgoing
mass flux as: 
 \begin{equation}
   \int_0^1\int_{v_z<0}\varphi_{out}(y,v)|v_z| \, dydv =\int_0^1\int_{v'_z>0}\varphi_{in}(y',v')|\Lambda_2 (y',v')| J\, dy'dv',  \label{chgtvar}
\end{equation}

Moreover, it can easily been obtained that the average mass flux is
zero, which reads
\begin{equation}
\int_0^1\int_{v_z<0}\varphi_{out}(y,v)|v_z| \, dydv =\int_0^1\int_{v'_z>0}\varphi_{in}(y',v')|v'_z| \, dy'dv' .  \label{fluxnul}
\end{equation}
\new{Indeed, integrating
equation~(\ref{LKSL-r-1-app}) with respect to $(y,z)$ on the cell $\{ y\in
[0,1], 0< z<\zeta_{\infty}(y)\}$, using the y-periodicity of
$\varphi$, and then taking into account that $\varphi$ is zero at
$z=\zeta_{\infty}(y)$, we get
 \begin{equation*}   \int_0^1 \varphi(y,0,v) v_z dy   =  - \frac{1}{2}\int 
\begin{pmatrix}  
\frac{1}{\beta_*}\partial_y{\cal V}_{\#} \\
\partial_z{\cal V}_{\#}
\end{pmatrix}
.\nabla_v\,  \varphi(y,0,v)\ dydz .  \end{equation*}
Then, we can integrate the previous relation with respect to $v$: the
right-hand side vanishes, and we get: 
 \begin{equation*}  
\int \int_0^1 \varphi(y,0,v) v_z \, dydv  =0,
\end{equation*}
which gives~(\ref{fluxnul}).
}

Now, we compare~(\ref{chgtvar}) and~(\ref{fluxnul}) to get
\begin{equation}
\int_0^1\int_{v'_z>0}\varphi_{in}(y',v')|\Lambda_2 (y',v')| J\, dy'dv'= \int_0^1\int_{v'_z>0}\varphi_{in}(y',v')|v'_z| \, dy'dv',
\end{equation}
which is true for every function $\varphi_{in}$. Consequently, we
deduce that the Jacobian $J$ satisfies
\begin{equation}
J= \frac{|v'_z|}{|\Lambda_2 (y',v')| }= \frac{|v'_z|}{|v_z|},
\end{equation}
which reads in the following more symmetric way
\begin{equation}
|v'_z|dy'dv' = |v_z|dydv.
\end{equation}

\section{Reciprocity property for the Maxwell like boundary
  condition~(\ref{roughmaxlike})}
\label{app:recip_maxwell}

\paragraph{Computation of $\sigma$.}

This parameter can be determined with the constraint of zero mass flux
of $F^0$ through the boundary $z=0$, that is to say
\begin{equation*}
  \int_{v_z<0}v_zF^0(t,x,0,v) \, dv 
+ \int_{v'_z>0}v'_zF^0(t,x,0,v') \, dv' =0.
\end{equation*}
\new{
Indeed, integrating
equation~(\ref{LKSL-r-1-bis}) with respect to $(y,z)$ on the cell $\{ y\in
[0,1], 0< z<\zeta_{\infty}(y)\}$, using the y-periodicity of
$\phi^0$, and then taking into account that $\phi^0$ is zero at
$z=\zeta_{\infty}(y)$, we get
 \begin{equation*}   \int_0^1 \phi^0 v_z dy   =  - \frac{1}{2}\int 
\begin{pmatrix}  
\frac{1}{\beta_*}\partial_y{\cal V}_{\#} \\
\partial_z{\cal V}_{\#}
\end{pmatrix}
.\nabla_v \phi^0\ dydz .  \end{equation*}
 But~(\ref{LKSL-r-2-bis}) and~(\ref{boltzmann-r-2-bis}) imply $\int_0^1 \phi^0 v_z dy =v_zF^0(t,x,0,v_x,v_z)$, so that after integration in $v$ we obtain
 \begin{equation}   \int F^0 v_z dv   = -  \int 
\begin{pmatrix}  
\frac{1}{\beta_*}\partial_y{\cal V}_{\#} \\
\partial_z{\cal V}_{\#}
\end{pmatrix}
.\left(\int \nabla_v \phi^0\  dv \right)\,  dydz =0,  \label{zeroflux-r-bis} \end{equation}
 which means that the mass flux of $F^0$  through the boundary $z=0$
vanishes.
}

Then using~(\ref{eq-routgoing}), we rapidly find
\begin{equation*}
\sigma(t,x)=-\frac{\int_{v'_z>0} 
\left( v'_z + \int_{v_z<0} v_zk_1(v'\rightarrow v)\, dv  
             \right) F^0(t,x,0,v')\, dv'}
{
\int_{v_z<0}v_z
\left(  1-\int_0^1 \exp(-r(y,-v))\, dy  \right) M(v) \, dv }.
\end{equation*}

Note that the integral with $k_1$ can be computed: by
using~(\ref{eq-k1}), the change of variables $(y,v)=\Lambda(y',w')$ and the
property of $r(y',v')$ mentioned above, we find
\begin{equation*}
\begin{split}
\int_{v_z<0} v_zk_1(v'\rightarrow v)\, dv 
& = -\int_{v_z<0}\int_0^1|v_z| \exp(-r(y,-v))\delta(v'+\Lambda_2(y,-v))\,
dydv \\ 
& = -\int_{w'_z>0}\int_0^1|w'_z| \exp(-r(y',w'))\delta(v'-w')\,
dy'dw' \\ 
& = -v'_z \int_0^1 \exp(-r(y',v')) \, dy' .
\end{split}
\end{equation*}

Consequently, the final form of $\sigma$ is:
\begin{equation}  \label{eq-sigma}
\sigma(t,x)=\frac{\int_{v'_z>0} v'_z
\left(1  -  \int_0^1 \exp(-r(y',v')) \, dy'
             \right) F^0(t,x,0,v')\, dv'}
{-
\int_{v_z<0}v_z
\left(  1-\int_0^1 \exp(-r(y,-v))\, dy  \right) M(v) \, dv}, 
\end{equation}
where the denominator is a constant denoted by $C$ in the following.

\paragraph{Scattering kernel for the boundary condition~(\ref{eq-routgoing}).}
Using~(\ref{eq-sigma}) in~(\ref{eq-routgoing}), we find
\begin{equation*}
\begin{split}
F^0(t,x,0,v)=
& \int_{v'_z>0}k_1(v'\rightarrow v)F^0(t,x,0,v')\, dv' \\
&  +  \left(1-\int_0^1\exp (-r(y,-v))\, dy \right)\, 
\frac{1}{C}\int_{v'_z>0} v'_z
\left(1  -  \int_0^1 \exp(-r(y',v')) \, dy'
             \right) F^0(t,x,0,v')\, dv'
M(v)  \\
& = \int_{v'_z>0}k^\#(v'\rightarrow v)F^0(t,x,0,v')\, dv',
\end{split}
\end{equation*}
with the scattering kernel $k^{\#}(v'\rightarrow v)=k_1(v'\rightarrow v)+k_2(v'\rightarrow
v)$, where
\begin{equation*}  
k_2(v'\rightarrow v) = 
\left(1-\int_0^1\exp (-r(y,-v))\, dy \right)\, 
\frac{1}{C} \,  v'_z
\left(1  -  \int_0^1 \exp(-r(y',v')) \, dy'
             \right) M(v).
\end{equation*}
This kernel can be written
\begin{equation}  \label{eq-defk2}
k_2(v'\rightarrow v) = 
\frac{1}{C} \psi(-v)\psi(v') |v'_z|M(v),
\end{equation}
where $\psi(w)=1  -  \int_0^1 \exp(-r(y,w)) \, dy$.

Consequently, the reciprocity of $k^\#$ can be deduced from the
reciprocity of the kernels $k_1$ and $k_2$.

\paragraph{Reciprocity of $k_1$.}

Using the definition of $k_1$~((\ref{eq-k1})), we have
\begin{equation}\label{eq-vk1M} 
  |v_z|k_1(v'\rightarrow v)
  M(v')=\int_0^1\exp(-r(y,-v))\delta(v'+\Lambda_2(y,-v))|v_z|M(v')\, dy.
\end{equation}
Then for a given $v$ and some test function $\theta$, we have
\begin{equation}\label{eq-recipk1_1} 
\begin{split}
\int_{v'_z>0} |v_z|k_1(v'\rightarrow v) M(v') \theta(v')\, dv'
& =
\int_{v'_z>0}\int_0^1\exp(-r(y,-v))\delta(v'+\Lambda_2(y,-v))|v_z|M(v')\theta(v')\,dydv'\\
& =
\int_0^1\exp(-r(y,-v))|v_z|M(-\Lambda_2(y,-v))\theta(-\Lambda_2(y,-v))\,dy\\
& =
\int_0^1\exp(-r(y,-v))\theta(-\Lambda_2(y,-v))\,dy\, |v_z|M(v),
\end{split}
\end{equation}
where we used~(\ref{eq-charnorme}). 

Moreover, we can use~(\ref{eq-vk1M}) to write
\begin{equation*}
  |v'_z|k_1(-v\rightarrow -v')
  M(v)=\int_0^1\exp(-r(y',v'))\delta(-v+\Lambda_2(y',v'))|v'_z|M(v)\, dy'.
\end{equation*}
Then, with the same $v$ and test function $\theta$ as above, we have
\begin{equation}\label{eq-recipk1_2} 
\begin{split}
\int_{v'_z>0} |v'_z|k_1(-v\rightarrow -v') M(v) \theta(v')\, dv'
& =
\int_{v'_z>0}\int_0^1\exp(-r(y',v'))\delta(-v+\Lambda_2(y',v'))|v'_z|M(v)\theta(v')\,dy'dv'\\
& =\int_{w_z<0}\int_0^1\exp(-r(y,-w))\delta(-v+w)|w_z|M(v)\theta(-\Lambda_2(y,-w))\,dydw\\
& =\int_0^1\exp(-r(y,-v))\theta(-\Lambda_2(y,-v))\,dy\, |v_z|M(v).
\end{split}
\end{equation}

Comparing~(\ref{eq-recipk1_1}) and~(\ref{eq-recipk1_2}), we find that
the two left-hand sides are equal for every $v$ and every test
function $\theta$, and then we have
\begin{equation*}
|v_z|k_1(v'\rightarrow v) M(v') =   |v'_z|k_1(-v\rightarrow -v') M(v), 
\end{equation*}
which is the reciprocity relation for $k_1$.

\paragraph{Reciprocity of $k_2$.}

This property is straightforward: using~(\ref{eq-defk2}), we have
\begin{equation*}
\begin{split}
|v_z|  k_2(v'\rightarrow v) M(v')& = 
\frac{1}{C} \psi(-v)\psi(v') |v_z||v'_z|M(v)M(v') \\
& = \frac{1}{C} \psi(v')\psi(-v) |v'_z||v_z|M(-v')M(-v) \\
& = |v'_z|  k_2(-v\rightarrow -v') M(-v),
\end{split}
\end{equation*}
since $M$ depends only on the norm of $v$.

The reciprocity of $k^\#$ follows, which completes the proof.

\end{document}